\documentclass[final,1p,times]{elsarticle}

\usepackage{amssymb}
 \usepackage{amsthm}
\usepackage{amscd}
\usepackage{amsmath}
\usepackage{amsfonts}
\usepackage{amssymb}
\usepackage{graphicx}

\numberwithin{equation}{section}

\newcommand{\ra}{\rightarrow}
\newcommand{\p}{\partial}
\newcommand{\f}{\frac}
\newcommand{\g}{\gamma}

\newcommand{\be}{\begin{equation}}
\renewcommand{\ra}{\rightarrow}
\newcommand{\ee}{\end{equation}}
\newcommand{\bea}{\begin{eqnarray}}
\newcommand{\eea}{\end{eqnarray}}
\newcommand{\bna}{\begin{eqnarray*}}
\newcommand{\ena}{\end{eqnarray*}}

\renewcommand{\le}{\left}
\newcommand{\ri}{\right}

\journal{$\ast\ast\ast$}

\begin{document}

\begin{frontmatter}

\title{Adams type inequalities and related elliptic partial differential equations in dimension four}

\author{Yunyan Yang}
\ead{yunyanyang@ruc.edu.cn}
\address{Department of Mathematics,
Renmin University of China, Beijing 100872, P. R. China}

\begin{abstract}
Motivated by Ruf-Sani's recent work, we prove an Adams type inequality and a singular
Adams type inequality in the whole four dimensional Euclidean space. As applications of those inequalities,
a class of elliptic partial differential equations are considered. Existence of nontrivial weak
solutions and multiplicity results are obtained via the
mountain-pass theorem and the Ekeland's variational principle. This
is a continuation of our previous work about singular
Trudinger-Moser type inequality.

\end{abstract}


\end{frontmatter}

\section{Introduction and main results}
Let $\Omega\subset\mathbb{R}^n$ be a smooth bounded domain. The
classical Trudinger-Moser inequality \cite{Moser,Pohozave,Trudinger}
says
 \be\label{T-M}\sup_{u\in W_0^{1,n}(\Omega),\,\|\nabla u\|_{L^n(\Omega)}\leq 1}
 \int_\Omega e^{\alpha |u|^{\f{n}{n-1}}}dx<\infty\ee
 for all $\alpha\leq \alpha_n=n\omega_{n-1}^{1/(n-1)}$, where
 $\omega_{n-1}$ is the area of the unit sphere in $\mathbb{R}^n$.
 $(\ref{T-M})$ is sharp in the sense that for any $\alpha>\alpha_n$,
 the integrals in (\ref{T-M}) are still finite,
 but the supremum of the integrals are infinite. (\ref{T-M}) plays
 an essential role in the study of the following partial
 differential equations
 \be\label{quasi}\le\{\begin{array}{lll}-{\rm div}(|\nabla u|^{n-2}\nabla u)=f(x,u)\,\,\,{\rm
 in}\,\,\,
 \Omega\\[1.5ex]
 u\in W_0^{1,n}(\Omega)\setminus\{0\},\end{array}\ri.\ee
where, roughly speaking, $f(x,u)$ behaves like $e^{|u|^{n/(n-1)}}$
as $|u|\ra\infty$. Problem (\ref{quasi}) and similar problems were
studied by many authors. Here we mention Atkinson-Peletier
\cite{AP}, Carleson-Chang \cite{CC}, Adimurthi et al.
\cite{Adi1}-\cite{Adi4}, de Figueiredo-Miyagaki-Ruf \cite{dMR},
Panda \cite{Panda2}, J. M. do \'O \cite{doo2}, de Figueiredo-do
\'O-Ruf \cite{ddR}, Silva-Soares \cite{SS}, Yang-Zhao
\cite{Yangzhao}, do \'O-Yang \cite{doYang}, Lam-Lu \cite{LL} and the references therein.

When $\Omega=\mathbb{R}^n$, the integrals in (\ref{T-M}) are infinite.
To get a Trudinger-Moser type inequality in this case, D. Cao \cite{Cao} proposed
the following: $\forall \alpha<4\pi$, $\forall M>0$,
\be\label{cao}\sup_{\int_{\mathbb{R}^2}|\nabla u|^2dx\leq
1,\,\int_{\mathbb{R}^2}u^2dx\leq M} \int_{\mathbb{R}^2}\le(e^{\alpha
u^2}-1\ri)dx<\infty,\,\,\ee which is equivalent to saying that for any
$\tau>0$ and $\alpha<4\pi$,
 \be\label{cao1}\sup_{\int_{\mathbb{R}^2}\le(|\nabla u|^2+\tau u^2\ri)dx\leq 1}
\int_{\mathbb{R}^2}\le(e^{\alpha u^2}-1\ri)dx<\infty.\ee (\ref{cao})
was independently generalized by Panda \cite{Panda} and J. M. do \'O
\cite{doo1} to $n$-dimensional case. Later Adachi and Tanaka \cite{AT}
gave another type of generalization. (\ref{cao}) and its high
dimensional generalizations were extensively used to study the
equation
$$-{\rm div}(|\nabla u|^{n-2}\nabla u)+V(x)|u|^{n-2}u=f(x,u)\quad{\rm in}\quad \mathbb{R}^n,$$
where  $f(x,u)$
behaves like $e^{\alpha |u|^{n/(n-1)}}$ as $|u|\ra\infty$. See for
examples \cite{Alves,Cao,doo1,doo3,Panda}.

Notice that (\ref{cao}) or (\ref{cao1}) is a subcritical
Trudinger-Moser type inequality in the whole Euclidean space. While
the critical inequality was obtained by B. Ruf \cite{Ruf} in dimension two and Li-Ruf
\cite{Li-Ruf} in general dimension. Using a simple variable substitution, Adimurthi-Sandeep
\cite{A-S} established a singular Trudinger-Moser inequality, which
is generalized to the whole $\mathbb{R}^n$ by Adimurthi-Yang
\cite{Adi-Yang}, namely \be\label{sing}
\sup_{\int_{\mathbb{R}^n}\le(|\nabla u|^n+\tau |u|^n\ri)dx\leq
1}\int_{\mathbb{R}^n}\f{e^{\alpha|u|^{n/(n-1)}}-\sum_{k=0}^{n-2}\f{1}{k!}|u|
^{kn/(n-1)}}{|x|^\beta}dx<\infty, \ee where $0\leq\beta<n$,
$\alpha/\alpha_n+\beta/n\leq 1$, $\tau$ is any fixed positive real
number. When $\beta=0$ and $\tau=1$, (\ref{sing}) is the standard
critical Trudinger-Moser type inequality \cite{Li-Ruf,Ruf}. In
\cite{Adi-Yang} we also employed (\ref{sing}) to obtain existence of
weak solutions to the equation
$$-{\rm div}(|\nabla u|^{n-2}\nabla u)+V(x)|u|^{n-2}u=\f{f(x,u)}{|x|^\beta}+\epsilon
h\quad{\rm in}\quad \mathbb{R}^n,$$
where $f(x,u)$ behaves like $e^{\alpha |u|^{n/(n-1)}}$ as
$|u|\ra\infty$, $\epsilon>0$, $h$ belongs to the dual space of
$W^{1,n}(\mathbb{R}^n)$. Similar problems
were also studied by J. M. do \'O and M. de Souza \cite{do-de} in the case $n=2$.\\

Our aim is to derive similar results to \cite{Adi-Yang} for
bi-Laplacian operator in dimension four. The essential tool will be
the Adams type inequality in the whole $\mathbb{R}^4$. Let
$\Omega\subset\mathbb{R}^4$ be a smooth bounded domain. As a
generalization of the Trudinger-Moser inequality, Adams inequality
\cite{Adams1} reads \be\label{adam}
 \sup_{u\in W_0^{2,2}(\Omega),\,\int_\Omega |\Delta u|^2dx\leq 1}\int_\Omega
 e^{32\pi^2u^2}dx<\infty.
\ee
This inequality was extended by Tasi \cite{Tas} (see also  Theorem 3.1 in \cite{R-S}), namely
\be\label{Tas}
 \sup_{u\in W^{2,2}(\Omega)\cap W_0^{1,2}(\Omega),\,\int_\Omega |\Delta u|^2dx\leq 1}\int_\Omega
 e^{32\pi^2u^2}dx<\infty.
\ee Also the integrals in (\ref{adam}) will be infinite when $\Omega$
is replaced by the whole $\mathbb{R}^4$. But B. Ruf and F. Sani
\cite{R-S}
 were able to establish the corresponding Adams type inequality in $\mathbb{R}^4$, say\\

 \noindent {\bf Theorem A} (Ruf-Sani). {\it There holds
 \be\label{Rs}\sup_{u\in W^{2,2}(\mathbb{R}^4),\,\int_{\mathbb{R}^4}(-\Delta u+u)^2dx\leq 1}
 \int_{\mathbb{R}^4}(e^{32\pi^2u^2}-1)dx<\infty.\ee
 Furthermore this inequality is sharp, i.e. if $32\pi^2$ is replaced by any $\alpha>32\pi^2$, then
 the supremum is infinite.}\\

 \noindent In fact they obtained more in \cite{R-S}, but here we focus on four dimensional case. Noticing that
 for all $u\in W^{2,2}(\mathbb{R}^4)$
 $$\int_{\mathbb{R}^4}(-\Delta u+u)^2dx=\int_{\mathbb{R}^4}(|\Delta u|^2+2|\nabla u|^2+u^2)dx,$$
 one can rewrite (\ref{Rs}) as
 \be\label{rs}
  \sup_{u\in W^{2,2}(\mathbb{R}^4),\,\int_{\mathbb{R}^4}(|\Delta u|^2+2|\nabla u|^2+u^2)dx\leq 1}
 \int_{\mathbb{R}^4}(e^{32\pi^2u^2}-1)dx<\infty.
 \ee

 One of our goals is the following:\\

 \noindent{\bf Theorem 1.1.} {\it Let $0\leq\beta<4$. Then for all $\alpha>0$
and $u\in W^{2,2}(\mathbb{R}^4)$, there holds
 \be\label{ss}\int_{\mathbb{R}^4}\f{e^{\alpha u^2}-1}{|x|^\beta}dx<\infty.\ee
 Furthermore, assume $\tau$ and $\sigma$ are two positive constants, we have for all $\alpha<
 32\pi^2\le(1-\f{\beta}{4}\ri)$,
 \be\label{criticalineq}\sup_{u\in W^{2,2}(\mathbb{R}^4),\,\int_{\mathbb{R}^4}
 (|\Delta u|^2+\tau|\nabla u|^2+\sigma u^2)dx\leq 1}
 \int_{\mathbb{R}^4}\f{e^{\alpha u^2}-1}{|x|^\beta}dx<\infty.\ee
  When
$\alpha>32\pi^2\le(1-\f{\beta}{4}\ri)$, the supremum is
  infinite.}\\

  We remark that the inequality (\ref{criticalineq}) in Theorem 1.1 is only subcritical case. How to establish it
  in the critical case $\alpha=32\pi^2\le(1-{\beta}/{4}\ri)$ is still open. In Section 2, we will show that
   Theorem 1.1 can be derived from the following:\\

  \noindent{\bf Theorem 1.2.} {\it For all $\alpha>0$
and $u\in W^{2,2}(\mathbb{R}^4)$, there holds
 \be\label{ss1}\int_{\mathbb{R}^4}\le(e^{\alpha u^2}-1\ri)dx<\infty.\ee
 For all constants $\tau>0$ and
$\sigma>0$, there holds
 \be\label{3.1}\sup_{u\in W^{2,2}(\mathbb{R}^4),\,\int_{\mathbb{R}^4}
 (|\Delta u|^2+\tau|\nabla u|^2+\sigma u^2)dx\leq 1}
 \int_{\mathbb{R}^4}\le(e^{32\pi^2 u^2}-1\ri)dx<\infty.\ee
 Furthermore this inequality is sharp, i.e. if $32\pi^2$ is replaced by any $\alpha>32\pi^2$, then
 the supremum is infinite.}\\

 \noindent Though the second part of Theorem 1.2 is similar to Theorem A, (\ref{ss1}) and (\ref{3.1}) are more suitable
 to use than (\ref{Rs}) or (\ref{rs}) when
 considering the related partial differential equations. This is also our next goal.
 Precisely Theorem 1.1 can be applied to study the existence of weak solutions to the following nonlinear equation
 \be\label{problem}
 \begin{array}{lll}\Delta^2 u-{\rm div}(a(x)\nabla u)+b(x)u=\f{f(x,u)}{|x|^\beta}+\epsilon h(x)\quad {\rm in}\quad
 \mathbb{R}^4.\\[1.5ex]
\end{array}
 \ee
 Here and throughout this paper we assume $0\leq\beta<4$, $a(x)$, $b(x)$ are two continuous functions satisfying\\

 \noindent $(A_1)$ there exist two positive constants $a_0$ and $b_0$ such that $a(x)\geq
 a_0$ and $b(x)\geq b_0$ for all $x\in \mathbb{R}^4$;\\[1.5ex]
 \noindent $(A_2)$  $\f{1}{b(x)}\in L^1(\mathbb{R}^4)$.
 \\

\noindent We also assume  the following growth condition on the
nonlinearity $f(x,s)$:\\

\noindent $(H_1)$ There exist constants $\alpha_0$, $b_1$, $b_2>0$
and $\gamma\geq 1$ such that for all $(x,s)\in
\mathbb{R}^4\times\mathbb{R}$,
$$|f(x,s)|\leq
b_1 |s|+b_2|s|^\gamma\le(e^{\alpha_0s^2}-1\ri).$$

\noindent $(H_2)$ There exists $\mu>2$ such that for all $x\in
\mathbb{R}^4$ and $s\not=0$,
$$0<\mu F(x,s)\equiv \mu\int_0^sf(x,t)dt\leq sf(x,s).$$

\noindent $(H_3)$ There exist constants $R_0$, $M_0>0$ such that for
all $x\in\mathbb{R}^4$ and $|s|\geq R_0$,
$$0<F(x,s)\leq M_0|f(x,s)|.$$

\noindent Define a function space
\be\label{1.15}E=\le\{u\in W^{2,2}(\mathbb{R}^4):\int_{\mathbb{R}^4}\le(|\Delta u|^2+a(x)|\nabla u|^2+b(x)u^2\ri)dx<\infty\ri\}.\ee
We say that $u\in E$ is a weak solution of problem (\ref{problem}) if
for all $\varphi\in E$ we have
$$\int_{\mathbb{R}^4}\le(\Delta u\Delta \varphi+a(x)\nabla u\nabla\varphi+b(x)u\varphi\ri)dx
=\int_{\mathbb{R}^4}\f{f(x,u)}{|x|^\beta}\varphi
dx+\epsilon\int_{\mathbb{R}^4}h\varphi dx,$$ where $h\in E^*$. Here
and in the sequel we denote the dual space of $E$ by $E^*$. For all
$u\in E$, we denote for simplicity the norm of $u$ by \be\label{En}
\|u\|_E=\le(\int_{\mathbb{R}^4}\le(|\Delta u|^2+a(x)|\nabla
u|^2+b(x)u^2\ri)dx\ri)^{1/2}. \ee

For $\beta: 0\leq\beta<4$, we define a singular eigenvalue by
\be\label{lamda}\lambda_\beta=\inf_{u\in E,\,u\not\equiv
0}\f{\|u\|_E^2} {\int_{\mathbb{R}^4}\f{u^2}{|x|^\beta}dx}.\ee
If $\beta=0$, then by $(A_1)$, obviously we have $\lambda_0\geq b_0>0$. If $0<\beta<4$,
 the continuous embedding of $W^{2,2}(\mathbb{R}^4)\hookrightarrow
 L^q(\mathbb{R}^4)$ $(\forall q\geq 2)$ together with the H\"older inequality implies
 \be\label{lam}\int_{\mathbb{R}^4}\f{u^2}{|x|^\beta}dx\leq\int_{|x|>1}u^2dx+\le(\int_{|x|\leq 1}|u|^{2t}dx\ri)^{1/t}
 \le(\int_{|x|\leq 1}\f{1}{|x|^{\beta t^\prime}}dx\ri)^{1/t^\prime}
 \leq C\|u\|_{W^{2,2}(\mathbb{R}^4)}^2,\ee
 where $1/t+1/t^\prime=1$, $0<\beta t^\prime<4$ and
  $\|u\|_{W^{2,2}(\mathbb{R}^4)}^2=\int_{\mathbb{R}^4}\le(|\nabla^2u|^2+|\nabla
  u|^2+u^2\ri)dx$.
  Standard elliptic estimates (see for example \cite{GT}, Chapter 9) imply that the above $W^{2,2}(\mathbb{R}^4)$-norm is
  equivalent to
  \be\label{norm}\|u\|_{W^{2,2}(\mathbb{R}^4)}=\le(\int_{\mathbb{R}^4}\le(|\Delta u|^2+|\nabla u|^2+u^2\ri)dx\ri)^{1/2}.\ee
  In the sequel, we use (\ref{norm}) as the norm of function in $W^{2,2}(\mathbb{R}^4)$.
  Combining (\ref{En}), (\ref{lam}), (\ref{norm}) and the assumption $(A_1)$, we have
  $\int_{\mathbb{R}^4}\f{u^2}{|x|^\beta}dx\leq C\|u\|_E^2$.
 Hence, by (\ref{lamda}), we conclude $\lambda_\beta>0$.\\

 When $\epsilon=0$, (\ref{problem}) becomes
 \be\label{h-zero}
 \Delta^2u-{\rm div}(a(x)\nabla u)+b(x)u=\f{f(x,u)}{|x|^\beta}.
 \ee
 Now we state an application of Theorem 1.1 as follows:\\

 \noindent{\bf Theorem 1.3.} {\it Assume that $a(x)$ and $b(x)$ are two continuous functions satisfying $(A_1)$ and $(A_2)$.
 $f:\mathbb{R}^4\times\mathbb{R}\ra\mathbb{R}$ is a continuous
 function and the hypothesis $(H_1)$, $(H_2)$ and $(H_3)$ hold. Furthermore we assume
 $$\leqno(H_4)\quad\quad\quad \limsup_{s\ra
0}\f{2|F(x,s)|}{s^2}<\lambda_\beta\,\,\,{\rm
uniformly\,\,with\,\,respect\,\,to}\,\,\,x\in \mathbb{R}^4;$$
 $$\leqno(H_5)\quad\quad\quad \liminf_{s\ra
+\infty}sf(x,s)e^{-\alpha_0s^2}=+\infty\,\,\,{\rm
uniformly\,\,with\,\,respect\,\,to}\,\,\,x\in \mathbb{R}^4.$$
 Then the equation (\ref{h-zero}) has a nontrivial mountain-pass type weak solution $u\in
 E$.}\\

 We remark that the result in Theorem 1.3 is stronger than Theorem 1.2 of \cite{Adi-Yang} in
 the case $\epsilon=0$. One reason is that $E$ is compactly embedded in $L^q(\mathbb{R}^4)$ for
 all $q\geq 1$ (see Lemma 3.6 below), but $E$ is compactly embedded in $L^q(\mathbb{R}^N)$ for
 all $q\geq N$ under the assumptions in \cite{Adi-Yang}. The other reason is that here we have the
 additional assumption $(H_5)$.\\

 When $\epsilon\not= 0$, we have the following:\\

 \noindent{\bf Theorem 1.4.} {\it Assume that $a(x)$ and $b(x)$ are two continuous functions satisfying $(A_1)$ and $(A_2)$.
 $f:\mathbb{R}^4\times\mathbb{R}\ra\mathbb{R}$ is a continuous
 function and the hypothesis $(H_1)$, $(H_2)$ and $(H_3)$ hold. Furthermore we assume $(H_4)$.
Then there exists $\epsilon_1>0$ such that if $0<
\epsilon<\epsilon_1$, then the problem (\ref{problem}) has a weak
solution of mountain-pass type.
 }\\

\noindent {\bf Theorem 1.5.} {\it  Assume that $a(x)$ and $b(x)$ are
two continuous functions satisfying $(A_1)$ and $(A_2)$.
 $f:\mathbb{R}^4\times\mathbb{R}\ra\mathbb{R}$ is a continuous
 function and the hypothesis $(H_1)$, $(H_2)$ and $(H_4)$ hold.
 Furthermore assume $h\not\equiv 0$. Then there exists
$\epsilon_2>0$ such that if $0<\epsilon<\epsilon_2$, then the
problem (\ref{problem}) has a weak
solution with negative energy.}\\

The most interesting question is that under what conditions the two solutions
obtained in Theorem 1.4 and Theorem 1.5 are distinct. Precisely we have the following:\\

\noindent {\bf Theorem 1.6.} {\it Assume that $a(x)$ and $b(x)$ are
two continuous functions satisfying $(A_1)$ and $(A_2)$.
 $f:\mathbb{R}^4\times\mathbb{R}\ra\mathbb{R}$ is a continuous
 function and the hypothesis $(H_1)$, $(H_2)$, $(H_3)$, $(H_4)$ and $(H_5)$ hold.
 Furthermore assume $h\not\equiv 0$. Then there exists
$\epsilon_3>0$ such that if $0<\epsilon<\epsilon_3$, then the problem
(\ref{problem}) has two distinct weak solutions.}\\

Before ending the introduction, we give an example of $f(x,s)$
satisfying $(H_1)-(H_5)$, say
\be\label{fun}f(x,s)=\psi(x)s(e^{\alpha_0s^2}-1),\ee where
$\alpha_0>0$ and $\psi(x)$ is a continuous function with
$0<c_1\leq\psi\leq c_2$ for constants $c_1$ and $c_2$. Obviously
 $(H_1)$ is satisfied. Integrating (\ref{fun}), we have
\be
\label{F}F(x,s)=\int_0^sf(x,t)dt=\f{1}{2\alpha_0}\psi(x)\le(e^{\alpha_0s^2}-1-\alpha_0s^2\ri).\ee
For $2<\mu\leq 4$, we have for $s\not=0$,
 \bna
 0<\mu
 F(x,s)=\f{\mu}{2\alpha_0}\psi(x)\sum_{k=2}^\infty\f{\alpha_0^ks^{2k}}{k!}
 \leq\f{\mu}{4}\psi(x)s^2\sum_{k=2}^\infty\f{\alpha_0^{k-1}s^{2(k-1)}}{(k-1)!}\leq
 sf(x,s).
 \ena
 Hence $(H_2)$ holds.
 It follows from (\ref{fun}) and (\ref{F}) that
 $0<F(x,s)\leq \f{1}{2\alpha_0}f(x,s)$ for $|s|\geq 1$. Thus
 $(H_3)$ is satisfied. By (\ref{F}), we have ${F(x,s)}/{s^2}\ra 0$ as $s\ra 0$.
 Hence $(H_4)$ holds. Finally $(H_5)$ follows from (\ref{fun}) immediately.\\

We organize this paper as follows: In Section 2, we prove an Adams
type inequality and a singular Adams type inequality
 in the whole $\mathbb{R}^4$ (Theorem 1.1 and Theorem 1.2). Applications of singular Adams inequality  (Theorems 1.3-1.6) will
be shown in Section 3.

\section{Adams type inequality in the whole $\mathbb{R}^4$}
In this section, we will prove Theorem 1.1 and Theorem 1.2. Let us
first prove Theorem 1.2 by using the density of
$C_0^\infty(\mathbb{R}^4)$ in $W^{2,2}(\mathbb{R}^4)$ and an
argument of Ruf-Sani \cite{R-S}. \\

\noindent{\it Proof of Theorem 1.2.} Firstly we prove (\ref{3.1}).
$\forall \tau>0$, $\sigma>0$, we denote
$c_0=\min\{\tau/2,\sqrt{\sigma}\}$. Let $u$ be a function belonging
to $W^{2,2}(\mathbb{R}^4)$ and satisfying
$$\int_{\mathbb{R}^4}(-\Delta u+c_0u)^2dx= 1,$$ or equivalently
 $$\int_{\mathbb{R}^4}(|\Delta u|^2+2c_0|\nabla u|^2+c_0^2u^2)dx= 1.$$
By the density of $C_0^\infty(\mathbb{R}^4)$ in $
W^{2,2}(\mathbb{R}^4)$, without loss of generality, we can find a
sequence of functions $u_k\in C_0^\infty(\mathbb{R}^4)$ such that
$u_k\ra u$ in $W^{2,2}(\mathbb{R}^4)$ as $k\ra\infty$ and
$\int_{\mathbb{R}^4}(-\Delta u_k+c_0u_k)^2dx= 1$. For otherwise we
can use
$$\widetilde{u}_k=\f{u_k}{\le(\int_{\mathbb{R}^4}(-\Delta
u_k+c_0u_k)^2dx\ri)^{1/2}}$$ instead of $u_k$. Now suppose ${\rm
supp}\, u_k\subset \mathbb{B}_{R_k}$ for any fixed $k$. Let
$f_k=-\Delta u_k+c_0u_k$. Consider the problem $$\le\{
 \begin{array}{lll}
  -\Delta v_k+c_0v_k=f_k^\sharp\quad{\rm in}\quad \mathbb{B}_{R_k}\\
  [1.5ex] v_k\in W_0^{1,2}(\mathbb{B}_{R_k}),
 \end{array}
\ri.$$ where $f_k^\sharp$ is the Schwarz decreasing rearrangement of
$f_k$ (see for example \cite{Polya}). By the property of rearrangement, we have
 \be\label{nm}\int_{\mathbb{B}_{R_k}}(-\Delta v_k+c_0v_k)^2dx=
 \int_{\mathbb{B}_{R_k}}(-\Delta u_k+c_0u_k)^2dx=1.\ee
 It follows from Trombetti-Vazquez
\cite{T-V} that $v_k$ is radially symmetric and \be\label{symm}
 \int_{\mathbb{B}_{R_k}}(e^{32\pi^2 {u_k}^2}-1)dx=\int_{\mathbb{B}_{R_k}}(e^{32\pi^2 {u_k^\sharp}^2}-1)dx\leq
 \int_{\mathbb{B}_{R_k}}(e^{32\pi^2 {v_k}^2}-1)dx.
\ee The radial lemma (\cite{K}, Lemma 1.1, Chapter 6) implies
\be\label{radial}
 |v_k(x)|\leq
 \f{1}{\sqrt{2}\pi}\f{1}{|x|^{{3}/{2}}}\|v_k\|_{W^{1,2}(\mathbb{R}^4)}.
\ee The equality (\ref{nm}) implies that
 \be\label{w12}\|v_k\|_{W^{1,2}(\mathbb{R}^4)}=\le(\int_{\mathbb{B}_{R_k}}(|\nabla v_k|^2+v_k^2)dx\ri)^{1/2}\leq
 \sqrt{\f{1}{2c_0}+\f{1}{c_0^2}}.\ee
  Choose $r_0=\le(\f{1}{2\pi^2}\le(\f{1}{2c_0}+\f{1}{c_0^2}\ri)\ri)^{1/3}$.
 If $R_k\leq r_0$, then (\ref{Tas}) and (\ref{nm}) imply
 \be\label{out}
 \int_{\mathbb{B}_{R_k}}(e^{32\pi^2v_k^2}-1)dx\leq C
\ee for some constant $C$ depending only on $c_0$. If $R_k>r_0$,
(\ref{radial}) implies that $|v_k(x)|\leq 1$ when $|x|\geq r_0$.
Thus we have by (\ref{nm}),
 \be\label{gtr}
 \int_{\mathbb{B}_{R_k}\setminus\mathbb{B}_{r_0}}(e^{32\pi^2v_k^2}-1)dx\leq
 \sum_{j=1}^\infty\f{(32\pi^2)^j}{j!}\int_{\mathbb{B}_{R_k}}v_k^2dx\leq
 \f{1}{c_0^2}\sum_{j=1}^\infty\f{(32\pi^2)^j}{j!}.
\ee
 On $\mathbb{B}_{r_0}$, we have for any $\epsilon>0$ by using the Young inequality
 $$v_k^2(x)\leq (1+\epsilon)(v_k(x)-v_k(r_0))^2+\le(1+\f{1}{\epsilon}\ri)v_k^2(r_0).$$
 Take $\epsilon$ such that $$\f{1}{1+\epsilon}=\int_{\mathbb{B}_{R_k}}|\Delta v_k|^2dx
 =1-2c_0\int_{\mathbb{B}_{R_k}}|\nabla v_k|^2dx-c_0^2\int_{\mathbb{B}_{R_k}}
 v_k^2dx.$$
 It follows that
 $$1+\f{1}{\epsilon}=\f{1}{2c_0\int_{\mathbb{B}_{R_k}}|\nabla v_k|^2dx+c_0^2\int_{\mathbb{B}_{R_k}}
 v_k^2dx}\leq \f{1}{\min\{2c_0,c_0^2\}\|v_k\|^2_{W^{1,2}}}.$$
 This together with (\ref{radial}) and (\ref{w12}) gives
 $$\le(1+\f{1}{\epsilon}\ri)v_k^2(r_0)\leq \f{1}{\min\{2c_0,c_0^2\}}\f{1}{2\pi^2r_0^3}=
 \f{1}{\min\{2c_0,c_0^2\}\le(\f{1}{2c_0}+\f{1}{c_0^2}\ri)}<1.$$
 Notice that $v_k(x)-v_k(r_0)\in W^{2,2}(\mathbb{B}_{r_0})\cap
 W_0^{1,2}(\mathbb{B}_{r_0})$ and $\int_{\mathbb{B}_{R_k}}|\Delta v_k|^2dx\geq
 \int_{\mathbb{B}_{r_0}}|\Delta (v_k-v_k(r_0))|^2dx$,
 we obtain by (\ref{Tas})
 $$\int_{\mathbb{B}_{r_0}}(e^{32\pi^2v_k^2}-1)dx\leq C$$
 for some constant $C$ depending only on $c_0$. This together with
 (\ref{symm}), (\ref{out}), (\ref{gtr}) and Fatou's Lemma implies that there exists a constant $C$ depending only on
 $c_0$ such that
 \be\label{fn}\int_{\mathbb{R}^4}(e^{32\pi^2u^2}-1)dx\leq \liminf_{k\ra\infty}
 \int_{\mathbb{B}_{R_k}}(e^{32\pi^2u_k^2}-1)dx\leq \liminf_{k\ra\infty}
 \int_{\mathbb{B}_{R_k}}(e^{32\pi^2v_k^2}-1)dx\leq C.\ee
 Notice that
 $$\int_{\mathbb{R}^4}
 (|\Delta u|^2+\tau|\nabla u|^2+\sigma u^2)dx\geq \int_{\mathbb{R}^4}
 (-\Delta u+c_0 u)^2dx.$$
 We obtain
 \bna
 \sup_{\int_{\mathbb{R}^4}
 (|\Delta u|^2+\tau|\nabla u|^2+\sigma u^2)dx\leq 1}
 \int_{\mathbb{R}^4}\le(e^{32\pi^2 u^2}-1\ri)dx&\leq&
 \sup_{\int_{\mathbb{R}^4}
 (-\Delta u+c_0 u)^2dx\leq 1}
 \int_{\mathbb{R}^4}\le(e^{32\pi^2 u^2}-1\ri)dx\\
 &=&\sup_{\int_{\mathbb{R}^4}
 (-\Delta u+c_0 u)^2dx= 1}
 \int_{\mathbb{R}^4}\le(e^{32\pi^2 u^2}-1\ri)dx.
 \ena
 This together with (\ref{fn}) implies (\ref{3.1}).

 Secondly, for $\alpha>32\pi^2$, we employ a sequence of functions $u_\epsilon$ constructed in Section 2 of
 \cite{Lu-Yang} (see also (33) in \cite{R-S}). Let $\widetilde{u_\epsilon}=u_\epsilon/
 \le(\int_{\mathbb{R}^4}(|\Delta u_\epsilon|^2+\tau|\nabla u_\epsilon|^2+\sigma u_\epsilon^2)dx\ri)^{1/2}$. A straightforward
 calculation shows that
 $$\sup_{u\in W^{2,2}(\mathbb{R}^4),\,\int_{\mathbb{R}^4}
 (|\Delta u|^2+\tau|\nabla u|^2+\sigma u^2)dx\leq 1}
 \int_{\mathbb{R}^4}\le(e^{32\pi^2 u^2}-1\ri)dx\geq\int_{\mathbb{R}^4}\le(e^{32\pi^2\widetilde{u}_\epsilon^2}-1\ri)dx\ra +\infty$$
 as $\epsilon\ra 0$. Hence $32\pi^2$ is the best constant for (\ref{3.1}).

 Thirdly we prove (\ref{ss1}). Let
$\alpha>0$ be a real number and $u$ be a function belonging to
$W^{2,2}(\mathbb{R}^4)$. By the density of
$C_0^\infty(\mathbb{R}^4)$ in $W^{2,2}(\mathbb{R}^4)$, there exists
some $u_0\in C_0^\infty(\mathbb{R}^4)$ such that
$$\|u-u_0\|_{W^{2,2}(\mathbb{R}^4)}<\f{1}{\sqrt{2\alpha}}.$$
Here we use (\ref{norm}) as the definition of $W^{2,2}(\mathbb{R}^4)$-norm. Thus
$$\int_{\mathbb{R}^4}(|\Delta(u-u_0)|^2+|\nabla (u-u_0)|^2+(u-u_0)^2)dx\leq \f{1}{2\alpha}.$$
Assume ${\rm supp}\,u_0\subset\mathbb{B}_{R}$ for some $R>0$ and
$|u_0|\leq M$ for some $M>0$. Using the inequality $(a+b)^2\leq
2a^2+2b^2$, we have
\bea{\nonumber}
 \int_{\mathbb{R}^4}\le(e^{\alpha u^2}-1\ri)dx&\leq&\int_{\mathbb{R}^4}\le(e^{2\alpha
 (u-u_0)^2+2\alpha u_0^2}-1\ri)dx\\{\nonumber}
 &\leq&e^{2\alpha M^2}\int_{\mathbb{R}^4}\le(e^{2\alpha (u-u_0)^2}-1\ri)dx+
 \int_{\mathbb{R}^4}\le(e^{2 \alpha u_0^2}-1\ri)dx\\\label{28}
 &\leq&e^{2\alpha M^2}\int_{\mathbb{R}^4}\le(e^{2\alpha
 (u-u_0)^2}-1\ri)dx+(e^{2\alpha M^2}-1)|\mathbb{B}_R|,
\eea where $|\mathbb{B}_R|$ denotes the volume of $\mathbb{B}_R$.
 By (\ref{3.1}) with $\tau=1$ and $\sigma=1$, we have
 $$\int_{\mathbb{R}^4}\le(e^{2\alpha (u-u_0)^2}-1\ri)dx\leq C$$
 for some universal constant $C$. Thus (\ref{ss1}) follows from (\ref{28}) immediately.
 $\hfill\Box$\\

 Now we use the H\"older inequality
 and Theorem 1.2 to prove Theorem 1.1. To do this, we need a technical lemma, namely\\

 \noindent{\bf Lemma 2.1.} {\it For all $p\geq 1$ and $t\geq 1$, we have $(t-1)^p\leq t^p-1$. In particular
  $\le(e^{s^2}-1\ri)^p\leq e^{ps^2}-1$ for all $s\in\mathbb{R}$ and $p\geq 1$.}\\

  \noindent{\it Proof.} For all $p\geq 1$ and $t\geq 1$, we set
  $$\varphi(t)=t^p-1-(t-1)^p.$$
 Since the derivative of $\varphi$ satisfies
 $$\f{d}{dt}\varphi(t)=pt^{p-1}-p(t-1)^{p-1}\geq 0,\,\,\forall t\geq 1,$$
 thus $\varphi(t)\geq 0$ for all $t\geq 1$ and the lemma follows immediately.
 $\hfill\Box$\\

\noindent{\it Proof of Theorem 1.1}: For any $\alpha>0$ and $u\in
W^{2,2}(\mathbb{R}^4)$, we have by using the H\"older inequality and
Lemma 2.1 \bea\int_{\mathbb{R}^4}\f{e^{\alpha
u^2}-1}{|x|^\beta}dx&=&
 \int_{|x|>1}\f{e^{\alpha u^2}-1}{|x|^\beta}dx+\int_{|x|\leq 1}\f{e^{\alpha
 u^2}-1}{|x|^\beta}dx{\nonumber}\\{\nonumber}
 &\leq& \int_{\mathbb{R}^4}(e^{\alpha u^2}-1)dx+\le(\int_{|x|\leq 1}\le(e^{\alpha u^2}-1\ri)^pdx\ri)^{1/p}
 \le(\int_{|x|\leq 1}\f{1}{|x|^{\beta q}}dx\ri)^{1/q}\\\label{2.71}
 &\leq&\int_{\mathbb{R}^4}(e^{\alpha u^2}-1)dx+C\le(\int_{\mathbb{R}^4}\le(e^{\alpha p u^2}-1\ri)dx\ri)^{1/p}\eea
 for some constant $C$ depending only on $q$ and $\beta$, where $q>1$ is a real number
 such that $\beta q<4$ and $1/p+1/q=1$.
 This together with (\ref{ss1}) implies (\ref{ss}).

Assume $\alpha<32\pi^2(1-\beta/4)$ and $u\in W^{2,2}(\mathbb{R}^2)$
satisfies
$$\int_{\mathbb{R}^4}\le(|\Delta u|^2+\tau|\nabla u|^2+\sigma u^2\ri)dx\leq 1.$$
Coming back to (\ref{2.71}), since $\beta q<4$ and $1/p+1/q=1$, one has
$$\alpha p< 32\pi^2\f{1-\beta/4}{1-1/q}.$$
We can further choose $q$ sufficiently close to $4/\beta$ such that $\alpha
p<32\pi^2$. Hence (\ref{criticalineq}) follows from (\ref{2.71}) and
(\ref{3.1}) immediately. $\hfill\Box$

 \section{Partial differential equations related to Adams type inequality in $\mathbb{R}^4$}
 In this section, we will use the mountain-pass theory to discuss the existence of solutions to the problem
 (\ref{problem}). Precisely we will prove Theorems 1.3-1.6. Firstly we construct the functional framework corresponding to
 (\ref{problem}). Secondly we analyze the geometry of the functional. Thirdly we use the mountain-pass theory to prove
 Theorem 1.3 and Theorem 1.4. Finally we use compactness analysis to prove Theorem 1.5 and Theorem 1.6.
 Throughout this section we
 assume that $f:\mathbb{R}^4\times\mathbb{R}\ra\mathbb{R}$ is a continuous
 function.\\

 \subsection{The functional\\}

  Now we use the notations of Section 1. For $u\in W^{2,2}(\mathbb{R}^4)$, we define a functional
 $$J_\epsilon(u)=\f{1}{2}\int_{\mathbb{R}^4}\le(|\Delta u|^2+a(x)|\nabla u|^2+b(x)u^2\ri)dx-
 \int_{\mathbb{R}^4}\f{F(x,u)}{|x|^\beta}dx-\epsilon\int_{\mathbb{R}^4}hudx,$$
 where $h\in E^*$, the dual space of $E$ (see (\ref{1.15})).
 When $\epsilon=0$, we write
 $$J(u)=\f{1}{2}\int_{\mathbb{R}^4}\le(|\Delta u|^2+a(x)|\nabla u|^2+b(x)u^2\ri)dx-
 \int_{\mathbb{R}^4}\f{F(x,u)}{|x|^\beta}dx.$$
 Here $F(x,s)=\int_0^sf(x,s)ds$.
 Since we assume $f(x,s)$, $a(x)$, $b(x)$ are all continuous functions and $(A_1)$, $(A_2)$, $(H_1)$ hold,
 it follows from Theorem 1.1 that
 $J_\epsilon$ or $J$ is well defined and
  \be\label{c1}J_\epsilon,\,J\in\mathcal{C}^1(E,\mathbb{R}).\ee
  Let us explain how to show (\ref{c1}). It suffices to show that if $u_j\ra u_\infty$ in $E$, then
  $J_\epsilon(u_j)\ra J_\epsilon (u_\infty)$ and $J_\epsilon^\prime(u_j)\ra J_\epsilon^\prime (u_\infty)$ in $E^*$
  as $j\ra\infty$. We point out a crucial fact:  for all $q\geq 1$, $E$ is embedded in $L^q(\mathbb{R}^4)$
  compactly and postpone its proof to Lemma 3.6 below. By $(H_1)$,
  \be\label{h1}
   |F(x,u_j)|\leq b_1u_j^2+b_2|u_j|^{\gamma+1}\le(e^{\alpha_0u_j^2}-1\ri).
  \ee
  Firstly, since $\|u_j\|_E$ is bounded and $E\hookrightarrow L^q(\mathbb{R}^4)$ is compact for all $q\geq  1$, we may assume
   $u_j\ra u_\infty$ in $L^q(\mathbb{R}^4)$ for all $q\geq 1$. An easy computation gives
  \be\label{hi}
   \lim_{j\ra\infty}\int_{\mathbb{R}^4}\f{|u_j|^q}{|x|^\beta}dx=\int_{\mathbb{R}^4}\f{|u_\infty|^q}{|x|^\beta}dx
   \quad{\rm for\,\,all}\quad q\geq 1.
  \ee
  Nextly we claim that
  \be\label{hia}\lim_{j\ra\infty}\int_{\mathbb{R}^4}\f{|u_j|^{\gamma+1}\le(e^{\alpha_0u_j^2}-1\ri)}{|x|^\beta}dx
  =\int_{\mathbb{R}^4}\f{|u_\infty|^{\gamma+1}\le(e^{\alpha_0u_\infty^2}-1\ri)}{|x|^\beta}dx.\ee
  For this purpose, we define a function $\varphi:\mathbb{R}^4\times[0,\infty)\ra \mathbb{R}$ by
  $$\varphi(x,s)=\f{s^{\gamma+1}\le(e^{\alpha_0s^2}-1\ri)}{|x|^\beta}.$$
  By the mean value theorem
  \bea{\nonumber}|\varphi(x,|u_j|)-\varphi(x,|u_\infty|)|&\leq& |{\p \varphi}/{\p s}(\xi)||u_j-u_\infty|\\
  \label{www}&\leq&(\eta(|u_j|)+\eta(|u_\infty|))\f{|u_j-u_\infty|}{|x|^\beta},
  \eea
  where $\xi$ lies between $|u_j(x)|$ and $|u_\infty(x)|$, $\eta:[0,\infty)\ra \mathbb{R}$ is a function defined by
  $$\eta(s)=\le((\gamma+1)s^\gamma+2\alpha_0 s^{\gamma+2}\ri)\le(e^{\alpha_0s^2}-1\ri)+2\alpha_0 s^{\gamma+2}.$$
  Using Lemma 2.1 and the inequalities $(a+b)^2\leq 2a^2+2b^2$, $ab-1\leq \f{a^p-1}{p}+\f{b^r-1}{r}$, where $a$, $b\geq 0$, $\f{1}{p}+\f{1}{r}=1$, we have
  \bna
  \int_{\mathbb{R}^4}\le(e^{\alpha_0u_j^2}-1\ri)^qdx&\leq&\int_{\mathbb{R}^4}\le(e^{2\alpha_0q(u_j-u_\infty)^2+2\alpha_0q u_\infty^2}-1\ri)dx
  \\&\leq& \f{1}{p}\int_{\mathbb{R}^4}\le(e^{2\alpha_0qp(u_j-u_\infty)^2}-1\ri)dx+\f{1}{r}
  \int_{\mathbb{R}^4}\le(e^{2\alpha_0qr u_\infty^2}-1\ri)dx.
  \ena
  Recalling that $\|u_j-u_\infty\|_E\ra 0$ as $j\ra\infty$ and applying Theorem 1.2, we can see that
  $$\sup_{j}\int_{\mathbb{R}^4}\le(e^{\alpha_0u_j^2}-1\ri)^qdx<\infty,\quad\forall q\geq 1.$$
  This together with the compact embedding $E\hookrightarrow L^q(\mathbb{R}^4)$ for all $q\geq 1$ implies that
  $\eta(|u_j|)+\eta(|u_\infty|)$ is bounded in $L^q(\mathbb{R}^4)$ for all $q\geq 1$. By (\ref{hi}) and (\ref{www}), the H\"older
  inequality leads to
  $$\lim_{j\ra\infty}\int_{\mathbb{R}^4}|\varphi(x,|u_j|)-\varphi(x,|u_\infty|)|dx=0.$$
  Hence $(\ref{hia})$ holds. In view of (\ref{h1}), we obtain by using (\ref{hi}), (\ref{hia}) and the generalized Lebesgue's dominated theorem
  $$\lim_{j\ra\infty}\int_{\mathbb{R}^4}\f{F(x,u_j)}{|x|^\beta}dx=\int_{\mathbb{R}^4}\f{F(x,u_\infty)}{|x|^\beta}dx.$$
  Therefore $J_\epsilon(u_j)\ra J_\epsilon(u_\infty)$ as $j\ra \infty$. In a similar way we can prove
  $J_\epsilon^\prime(u_j)\ra J_\epsilon^\prime(u_\infty)$ in $E^*$ as $j\ra\infty$. Hence (\ref{c1}) holds.

 It is easy to see that
 the critical point $u_\epsilon$ of $J_\epsilon$ is a weak solution
 to (\ref{problem}) and the critical point $u$ of $J$ is a weak solution
 to (\ref{h-zero}). Thus, to find weak solutions to (\ref{problem}) or (\ref{h-zero}), it suffices to
 find critical points of $J_\epsilon$ or $J$ in the function space $E$.
  \subsection{The geometry of the functional\\}

In this subsection, we describe the geometry of the functional
$J_\epsilon$.
\\

  \noindent{\bf Lemma 3.1.} {\it Assume that $(H_2)$ and $(H_3)$ are satisfied. Then $J_\epsilon(tu)\ra -\infty$ as
  $t\ra +\infty$, for all compactly supported $u\in W^{2,2}(\mathbb{R}^4)\setminus \{0\}$.}\\

  \noindent{\it Proof.} Assume $u$ is supported in a bounded domain
  $\Omega$. Since $f(x,s)$ is continuous, in view of $(H_2)$, there
  exists constants $c_1$, $c_2>0$ such that $F(x,s)\geq
  c_1|s|^\mu-c_2$ for all $(x,s)\in\overline{\Omega}\times\mathbb{R}$. It then follows that
  \bna J_\epsilon(t u)&=&\f{t^2}{2}\int_{\Omega}\le(|\Delta u|^2+a(x)|\nabla u|^2+b(x)u^2\ri)dx-
 \int_{\Omega}\f{F(x,t u)}{|x|^\beta}dx-\epsilon\, t \int_{\Omega}hudx\\
 &\leq&\f{t^2}{2}\int_{\Omega}\le(|\Delta u|^2+a(x)|\nabla
 u|^2+b(x)u^2\ri)dx-c_1t^\mu \int_\Omega\f{|u|^\mu}{|x|^\beta}dx+O(t).
 \ena
 Then the lemma holds since $\mu>2$.$\hfill\Box$\\

  \noindent {\bf Lemma 3.2.} {\it Assume that $(A_1)$,
 $(H_1)$ and $(H_4)$ hold. Then there exists $\epsilon_1>0$
 such that for any $\epsilon:$ $0<\epsilon<\epsilon_1$, there
  exist $r_\epsilon>0$, $\vartheta_\epsilon>0$ such that
   $J_{\epsilon}(u)\geq\vartheta_\epsilon$ for all $u$ with $\|u\|_E=r_\epsilon$.
   Furthermore $r_\epsilon$ can be chosen such that $r_\epsilon\ra 0$ as $\epsilon\ra
   0$. When $\epsilon=0$, there
  exist $r_0>0$, such that if $r\leq r_0$, then there exists  $\vartheta>0$ depending only on $r$ such that
   $J(u)\geq\vartheta$ for all $u$ with $\|u\|_E= r$.}\\

   \noindent {\it Proof.} By $(H_4)$, there exist $\tau$, $\delta>0$ such
   that if $|s|\leq \delta$, then
   $$|F(x,s)|\leq \f{\lambda_\beta-\tau}{2}|s|^2$$
   for all $x\in\mathbb{R}^4$. By $(H_1)$, there holds for $|s|\geq \delta$
   \bea{\nonumber}
    |F(x,s)|&\leq&\int_0^{|s|}\le\{b_1t+b_2t^\gamma(e^{\alpha_0t^2}-1)\ri\}dt\\{\nonumber}
    &\leq&\f{b_1}{2}s^2+{b_2}|s|^{\gamma+1}(e^{\alpha_0s^2}-1)\\
    &\leq& C|s|^q(e^{\alpha_0 s^2}-1){\nonumber}
   \eea
   for any $q>\gamma+1\geq 2$, where $C$ is a constant depending only on $b_1$, $b_2$, $q$
   and $\delta$.
   Combining the above two inequalities, we obtain for all $s\in\mathbb{R}$
   \be\label{Cap-f}
    |F(x,s)|\leq \f{\lambda_\beta-\tau}{2}|s|^2+C|s|^q(e^{\alpha_0
    s^2}-1).
   \ee
   Recall that  $\|u\|_E$ and $\lambda_\beta$ are defined by (\ref{En}) and (\ref{lamda}) respectively. It follows from
   (\ref{Cap-f}) that
   \bea\label{121}
    J_\epsilon(u)&=&\f{1}{2}\|u\|_E^2-\int_{\mathbb{R}^4}\f{F(x,u)}{|x|^\beta}dx-\epsilon
    \int_{\mathbb{R}^4}hudx{\nonumber}\\
    &\geq&\f{1}{2}\|u\|_E^2-
    \f{\lambda_\beta-\tau}{2\lambda_\beta}\|u\|_E^2-C
    \int_{\mathbb{R}^4}\f{e^{\alpha_0
    u^2}-1}{|x|^\beta}|u|^qdx-\epsilon\int_{\mathbb{R}^4}hudx{\nonumber}\\
    &\geq&\f{\tau}{2\lambda_\beta}\|u\|_E^2-C
    \int_{\mathbb{R}^4}\f{e^{\alpha_0
    u^2}-1}{|x|^\beta}|u|^qdx-\epsilon\|h\|_{E^*}\|u\|_E,
   \eea
   where $$\|h\|_{E^*}=\sup_{\|\varphi\|_E=1}\le|\int_{\mathbb{R}^4}h\varphi dx\ri|.$$
   Using the H\"older inequality, Lemma 2.1 and the continuous embedding $E\hookrightarrow L^q(\mathbb{R}^4)$,
   we have
   \bea\label{111}
    \int_{\mathbb{R}^4}\f{e^{\alpha_0
    u^2}-1}{|x|^\beta}|u|^qdx&\leq&\le(\int_{\mathbb{R}^4}\f{(e^{\alpha_0
    u^2}-1)^{r\,^\prime}}{|x|^{\beta
    r\,^\prime}}dx\ri)^{1/r\,^\prime}\le(\int_{\mathbb{R}^4}|u|^{qr}dx\ri)^{1/r}{\nonumber}\\
    &\leq&C\le(\int_{\mathbb{R}^4}\f{e^{\alpha_0 r\,^\prime
    u^2}-1}{|x|^{\beta
    r\,^\prime}}dx\ri)^{1/r\,^\prime}\|u\|_E^q,
   \eea
   where $1/r+1/r\,^\prime=1$, $0\leq\beta r\,^\prime<4$ and $C$ is a constant such that
   $\|u\|_{L^{qr}(\mathbb{R}^4)}\leq C^{1/q}\|u\|_E$.
   Here and in the sequel we often denote various constants by the same $C$.
   By $(A_1)$,
   $$\int_{\mathbb{R}^4}\le(|\Delta u|^2+a_0|\nabla u|^2+b_0u^2\ri)dx\leq
   \|u\|_E^2.$$
   Theorem 1.1 implies that if
   \be\label{str}\|u\|_E^2<\f{16\pi^2}{\alpha_0 r\,^\prime}\le(1-\f{\beta
   r\,^\prime}{4}\ri),\ee then
   $$\int_{\mathbb{R}^4}\f{e^{\alpha_0 r\,^\prime
    u^2}-1}{|x|^{\beta
    r\,^\prime}}dx\leq C$$
    for some constant $C$ depending only on $\alpha_0$, $\beta$ and
    $r\,^\prime$. This together with (\ref{111}) gives
    \be\label{bd}\int_{\mathbb{R}^4}\f{e^{\alpha_0
    u^2}-1}{|x|^\beta}|u|^qdx\leq C\|u\|_E^q,\ee
    provided that $u$ satisfies (\ref{str}). Hence, assuming (\ref{str}), we obtain by combining (\ref{121}) and (\ref{bd})
    \be\label{lower}J_\epsilon(u)\geq \|u\|_E\le(\f{\tau}{2\lambda_\beta}\|u\|_E-C\|u\|_E^{q-1}-\epsilon\|h\|_{E^*}\ri).\ee
    Since $\tau>0$, there holds for sufficiently small $r>0$,
    $$\f{\tau}{2\lambda_\beta}r-Cr^{q-1}\geq \f{\tau}{4\lambda_\beta}r.$$
    If $h\not\equiv 0$, for sufficiently small $\epsilon>0$, we may take $r_\epsilon$ and $\vartheta_\epsilon$ such that
    $$\f{\tau}{4\lambda_\beta}r_\epsilon=2\epsilon\|h\|_{E^*},\quad \vartheta_\epsilon=\epsilon r_\epsilon\|h\|_{E^*}.$$
    This implies $J_\epsilon(u)\geq \vartheta_\epsilon$ for all $u$
    with $\|u\|_E=r_\epsilon$ and $r_\epsilon\ra 0$ as $\epsilon\ra 0$. If $\epsilon=0$,
    (\ref{lower}) implies
    $$J(u)\geq \|u\|_E\le(\f{\tau}{2\lambda_\beta}\|u\|_E-C\|u\|_E^{q-1}\ri).$$
    Hence there exists some $r_0>0$ such that if $r\leq r_0$, then
    $J(u)\geq \f{\tau r^2}{4\lambda_\beta}$ for all $u$ with
    $\|u\|_E=r$.
     $\hfill\Box$\\

    \noindent{\bf Lemma 3.3.} {\it Assume $h\not\equiv 0$, $(A_1)$ and $(H_1)$ hold.
    Then there exist $\tau>0$ and $v\in E$ with $\|v\|_E=1$ such that
    $J_{\epsilon}(tv)<0$ for all $t$: $0<t<\tau$. Particularly
    $\inf_{\|u\|_E\leq \tau}J_{\epsilon}(u)<0$.}\\

    \noindent{\it Proof}. For any fixed $h\in E^*$, one can view $h$ as a linear functional defined on
    $E$ by
    $$\int_{\mathbb{R}^4}hudx,\quad \forall u\in E.$$
    By $(A_1)$, $E$ is a Hilbert space under the inner product
    $$\langle u,v\rangle=\int_{\mathbb{R}^4}
    \le(\Delta u\Delta v+a(x)\nabla u\nabla v+b(x)uv\ri)dx.$$
     By the Riesz representation theorem,
    $$
 \Delta^2 u-{\rm div}(a(x)\nabla u)+b(x)u= \epsilon h\quad {\rm in}\quad
 \mathbb{R}^4
 $$
 has a unique weak solution $u\in E$. If $h\not\equiv 0$, then for any fixed $\epsilon>0$, we have $u\not\equiv 0$ and
 $$\epsilon \int_{\mathbb{R}^4}hudx=\|u\|_E^2>0.$$
 A simple calculation shows
 \be\label{der}\f{d}{dt}J_\epsilon(tu)=t\|u\|_E^2-\int_{\mathbb{R}^4}
 \f{f(x,tu)}{|x|^\beta}udx-\epsilon\int_{\mathbb{R}^4}hudx.\ee
 By $(H_1)$, we have
 \bea\le|\int_{\mathbb{R}^4}\f{f(x,tu)}{|x|^\beta}udx\ri|&\leq& b_1|t|\int_{\mathbb{R}^4}\f{u^2}{|x|^\beta}dx
 +b_2|t|^\gamma\int_{\mathbb{R}^4}\f{e^{\alpha_0
 t^2u^2}-1}{|x|^\beta}|u|^{1+\gamma}dx{\nonumber}\\\label{3.10}
 &\leq& \f{b_1|t|}{\lambda_\beta}\|u\|_E^2+b_2|t|^\gamma\int_{\mathbb{R}^4}\f{e^{\alpha_0 t^2u^2}-1}{|x|^\beta}
 |u|^{1+\gamma}dx.\eea
 Using the same argument we prove (\ref{bd}), there exists some $t_0>0$ such that if $|t|<t_0$, then
 $$\int_{\mathbb{R}^4}\f{e^{\alpha_0 t^2u^2}-1}{|x|^\beta}|u|^{\gamma+1}dx\leq C$$
 for some constant $C$ depending only on $t_0$, $\alpha_0$ and $\beta$. It then
 follows from (\ref{3.10}) that
 $$\lim_{t\ra 0}\int_{\mathbb{R}^4}\f{f(x,tu)}{|x|^\beta}udx=0.$$
 This together with (\ref{der}) implies that
 there exists some $\delta>0$ such that
 $$\f{d}{dt}J_\epsilon(tu)<0,$$
 provided that $0<t<\delta$. Notice that $J_\epsilon(0)=0$, we have
 $J_\epsilon(tu)<0$ for all $0<t<\delta$. $\hfill\Box$

 \subsection{Min-max level\\}

 In this subsection, we estimate the min-max level of $J_\epsilon$ or $J$. To do this, we
 define a sequence of functions $\widetilde{\phi}_n$ by
 $$
 \widetilde{\phi}_n(x)=\le\{
 \begin{array}{lll}
  \sqrt{\f{\log n}{8\pi^2}}-\f{n^2}{\sqrt{32\pi^2\log
  n}}{|x|^2}+\f{1}{\sqrt{32\pi^2\log n}}&{\rm for} &|x|\leq
  {1}/{n}\\[1.5ex] \f{1}{\sqrt{8\pi^2\log
  n}}\log\f{1}{|x|}&{\rm for}&{1}/{n}<|x|\leq 1,\\[1.5ex]
  \zeta_n(x)&{\rm for}&|x|> 1
 \end{array}
 \ri.
 $$
 where $\zeta_n\in C_0^\infty(\mathbb{B}_{2}(0))$,
 $\zeta_n\mid_{\p\mathbb{B}_1(0)}=\zeta_n\mid_{\p\mathbb{B}_{2}(0)}=0$,
 $\f{\p\zeta_n}{\p\nu}\mid_{\p\mathbb{B}_1(0)}=\f{1}{\sqrt{8\pi^2\log
 n}}$, $\f{\p\zeta_n}{\p\nu}\mid_{\p\mathbb{B}_{2}(0)}=0$, and
 $\zeta_n$, $|\nabla\zeta_n|$, $\Delta\zeta_n$ are all $O(1/\sqrt{\log
 n})$. One can check that $\widetilde{\phi}_n\in W_0^{2,2}(\mathbb{B}_{2}(0))\subset
 W^{2,2}(\mathbb{R}^4)$. Straightforward calculations show that
 $$\|\widetilde{\phi}_n\|_2^2=O(1/{\log n}),\,\,\|\nabla\widetilde{\phi}_n\|_2^2=
 O(1/{\log n}),\,\,\|\Delta\widetilde{\phi}_n\|_2^2=1+O(1/{\log n})$$
 and thus
 $$\|\widetilde{\phi}_n\|_E^2=1+O(1/{\log n}).$$
 Set
 $$\phi_n(x)=\f{\widetilde{\phi}_n(x)}{\|\widetilde{\phi}_n\|_E}$$
 so that $\|\phi_n\|_E=1$. It is not difficult to see that
 \be\label{local}
  \phi_n^2(x)\geq {\f{\log n}{8\pi^2}}+O(1)\quad{\rm for}\quad
  |x|\leq 1/n.
 \ee

 \noindent{\bf Lemma 3.4.} {\it Assume $(H_2)$ and $(H_3)$. There exists a sufficiently large $\nu_0>0$ such that if
 $$\liminf_{s\ra+\infty}sf(x,s)e^{-\alpha_0s^2}>\nu_0$$
 uniformly with respect to $x\in\mathbb{R}^4$,
 then there exists some $n\in\mathbb{N}$ such that
 $$\max_{t\geq 0}\le(\f{t^2}{2}-\int_{\mathbb{R}^4}\f{F(x,t\phi_n)}{|x|^\beta}dx\ri)<\f{16\pi^2}{\alpha_0}
 \le(1-\f{\beta}{4}\ri).$$}

 \noindent{\it Proof.} Suppose by contradiction that for any large $\nu>0$
 \be\label{hyp1}\liminf_{s\ra+\infty}sf(x,s)e^{-\alpha_0s^2}>\nu\ee
 uniformly with respect to $x\in\mathbb{R}^4$,
 but for all $n\geq 2$
 \be\label{contradict}\max_{t\geq
 0}\le(\f{t^2}{2}-\int_{\mathbb{R}^4}\f{F(x,t\phi_n)}{|x|^\beta}dx\ri)\geq\f{16\pi^2}{\alpha_0}\le(1-\f{\beta}{4}\ri).\ee
 By $(H_2)$, we have $F(x,t\phi_n)>0$ and $f(x,t\phi_n)>0$ when
 $t\phi_n(x)>0$. In addition, $F(x,t\phi_n)\geq 0$ for all $x\in\mathbb{R}^4$. Furthermore, by $(H_3)$,
 there exists some constants $C_1$, $R_0$ and $M_0$ such that
 $F(s)\geq C_1 e^{s/M_0}$ for $s\geq R_0$. Hence we have
 \be\label{F-1}\int_{\mathbb{R}^4}\f{F(x,t\phi_n)}{|x|^\beta}dx\geq C_1\int_{t\phi_n\geq R_0}
 \f{e^{t\phi_n/M_0}}{|x|^\beta}dx.\ee
 For each fixed $n$, one can choose sufficiently large $T_n$ such that if $t\geq T_n$,
 then $t\phi_n\geq R_0$ on $\mathbb{B}_{R_0}(0)$. Thus for $t\geq
 T_n$ we have
 \be\label{F-2}\int_{t\phi_n\geq R_0}\f{e^{t\phi_n/M_0}}{|x|^\beta}dx\geq \f{t^3}{6M_0^3}\int_{t\phi_n\geq R_0}
 \f{\phi_n^3}{|x|^\beta}dx\geq \f{t^3}{6M_0^3}\int_{|x|\leq r/n}\f{\phi_n^3}{|x|^\beta}dx.\ee
 Combining (\ref{F-1}) and (\ref{F-2}) we have
\be\label{minus}\lim_{t\ra
+\infty}\le(\f{t^2}{2}-\int_{\mathbb{R}^4}\f{F(x,t\phi_n)}{|x|^\beta}dx\ri)=-\infty.\ee
 Since $F(x,0)=0$, it follows from (\ref{contradict}) and
 (\ref{minus}) that there exists $t_n>0$ such that
 \be\label{cr}\f{t_n^2}{2}-\int_{\mathbb{R}^4}\f{F(x,t_n\phi_n)}{|x|^\beta}dx=
 \max_{t\geq 0}\le(\f{t^2}{2}-\int_{\mathbb{R}^4}\f{F(x,t\phi_n)}{|x|^\beta}dx\ri).\ee
 Clearly there holds at $t=t_n$
 $$\f{d}{dt}\le(\f{t^2}{2}-\int_{\mathbb{R}^4}\f{F(x,t\phi_n)}{|x|^\beta}dx\ri)=0.$$
 Hence
 $$
 t_n^2=\int_{\mathbb{R}^4}\f{t_n\phi_n f(x,t_n\phi_n)}{|x|^\beta}dx.
 $$

 Now we prove that $\{t_n\}$ is bounded.
 Suppose not. By (\ref{local}), (\ref{hyp1}), (\ref{contradict}) and (\ref{cr}), there holds for
 sufficiently large $n$
 \be\label{t-bd}
 t_n^2\geq \f{\nu}{2}\int_{|x|\leq
 1/n}\f{e^{\alpha_0(t_n\phi_n)^2}}{|x|^\beta}dx\geq \f{\nu}{2}\f{\omega_3}{4-\beta}
 \f{1}{n^{4-\beta}}e^{\alpha_0t_n^2\le(\f{\log
 n}{8\pi^2}+O(1)\ri)},
 \ee
 where $\omega_3$ is the area of the unit sphere in $\mathbb{R}^4$.
 It follows that
 $$1\geq \f{\nu}{2}\f{\omega_3}{4-\beta}n^{\alpha_0t_n^2\le(\f{1}{8\pi^2}+o(1)\ri)+\beta-4-
 \f{2}{\log n}\log t_n}.$$
 Letting $n\ra\infty$, we get a contradiction because the term on the right hand tends to $+\infty$.

 Next we prove that $t_n^2\ra \f{32\pi^2}{\alpha_0}\le(1-\f{\beta}{4}\ri)$ as $n\ra \infty$.
 By $(H_2)$, $F(x,t_n\phi_n)\geq 0$. It then follows from (\ref{contradict}) and
 (\ref{cr}) that
 \be\label{ge}{t_n^2}\geq \f{32\pi^2}{\alpha_0}\le(1-\f{\beta}{4}\ri).\ee
 Suppose that $\lim_{n\ra\infty}t_n^2>\f{32\pi^2}{\alpha_0}\le(1-\f{\beta}{4}\ri)$. By (\ref{t-bd})
 and $\{t_n\}$ is bounded we get a contradiction. Thus we have $t_n^2\ra \f{32\pi^2}{\alpha_0}\le(1-\f{\beta}{4}\ri)$ as $n\ra \infty$.

 By (\ref{hyp1}), there exists some $s_0>0$ such that $sf(s)e^{-\alpha_0s^2}\geq
 {\nu}/{2}$ for all $s\geq s_0$. Since (\ref{t-bd}) holds for sufficiently
 large $n$, we have by combining (\ref{t-bd}) and (\ref{ge}) that
 $$\nu\leq \f{8}{\omega_3}t_n^2e^{-\alpha_0t_n^2O(1)}.$$
 Letting $n\ra\infty$, we get $\nu\leq C$ for some constant $C$ depending only on $\alpha_0$. This contradicts the arbitrary of $\nu$
 and completes the proof of the lemma.
 $\hfill\Box$\\

 \subsection{Palais-Smale sequence\\}

  In this subsection, we want to show that the weak limit of a Palais-Smale sequence for $J_\epsilon$ is a weak solution
  of $(\ref{problem})$. To this end, we need the following convergence
  result, which is a generalization of Lemma 2.1 in \cite{dMR} and Lemma 3.7 in
  \cite{do-de}.\\

  \noindent{\bf Lemma 3.5.} {\it Let $f:\mathbb{R}^N\times\mathbb{R}\ra \mathbb{R}\,\,(N\geq 1)$ be a measurable
  function. Assume for any $\eta>0$, there exists some constant $c$ depending only on $\eta$ such that
  $|f(x,s)|\leq c|s|$ for all $(x,s)\in
  \mathbb{R}^N\times[-\eta,\eta]$. Let $\phi:\mathbb{R}^N\ra
  \mathbb{R}$ be a nonnegative measurable function,
   $(u_n)$ be a sequence of functions  with $u_n\ra u$ in $\mathbb{R}^N$ almost everywhere,
   $\phi u_n\ra \phi u$ strongly in
  $L^1(\mathbb{R}^N)$, $\phi f(x,u)\in L^1(\mathbb{R}^N)$
  and \be\label{jie}\int_{\mathbb{R}^N}\phi|f(x,u_n)u_n|dx\leq
  C\ee for all $n$. Then, up to a subsequence, we have $\phi f(x,u_n)\ra \phi f(x,u)$ strongly
  in $L^1(\mathbb{R}^N)$.}\\

  \noindent{\it Proof.} Since $u,\,\phi f(x,u)\in L^1(\mathbb{R^N})$, we have
  $$\lim_{\eta\ra+\infty}\int_{|u|\geq\eta}\phi|f(x,u)|dx=0.$$
  Let $C$ be the constant in (\ref{jie}). Given any $\epsilon>0$, one can select some $M>{C}/{\epsilon}$ such that
  \be\label{abs}
  \int_{|u|\geq M}\phi|f(x,u)|dx<\epsilon.
  \ee
  It follows from (\ref{jie}) that
  \be\label{abe}\int_{|u_n|\geq M}\phi|f(x,u_n)|dx\leq \f{1}{M}\int_{|u_n|\geq
  M}\phi|f(x,u_n)u_n|dx<\epsilon.\ee
  For all $x\in\{x\in\mathbb{R}^N:|u_n|<M\}$, by our assumption, there exists a constant $C_1$ depending only on
  $M$ such that $|f(x,u_n(x))|\leq
  C_1|u_n(x)|$. Notice that $\phi u_n\ra\phi u$ strongly in $L^1(\mathbb{R}^N)$ and $u_n\ra u$ almost everywhere
  in $\mathbb{R}^N$.
  By a generalized Lebesgue's dominated convergence theorem, up to a subsequence we obtain
  \be\label{domin}\lim_{n\ra
  \infty}\int_{|u_n|<M}\phi|f(x,u_n)|dx=\int_{|u|<M}\phi|f(x,u)|dx.\ee
  Combining (\ref{abs}), (\ref{abe}) and (\ref{domin}), we can find
  some $K>0$ such that when $n>K$,
  $$\le|\int_{\mathbb{R}^N}\phi|f(x,u_n)|dx-\int_{\mathbb{R}^N}\phi|f(x,u)|dx\ri|<3\epsilon.$$
  Hence $\|\phi f(x,u_n)\|_{L^1(\mathbb{R}^N)}\ra
  \|\phi f(x,u)\|_{L^1(\mathbb{R}^N)}$ as $n\ra \infty$. Since  $\phi u_n\ra
  \phi u$ in $\mathbb{R}^N$ almost everywhere, we get the desired result.
  $\hfill\Box$\\

  \noindent{\bf Lemma 3.6.} {\it Assume $(A_1)$ and $(A_2)$. Then we have the compact embedding
   $$E\hookrightarrow L^q(\mathbb{R}^4)\quad{\rm for \,\,all}\quad q\geq 1.$$
  }
  \noindent{\it Proof.} By $(A_1)$, the Sobolev embedding
  theorem implies the following continuous embedding
  $$E\hookrightarrow W^{2,2}(\mathbb{R}^4)\hookrightarrow L^q(\mathbb{R}^4)\quad
  {\rm for\,\,all}\quad 2\leq q<\infty.$$
  It follows from the H\"older inequality and $(A_2)$ that
  $$\int_{\mathbb{R}^4}|u|dx\leq\le(\int_{\mathbb{R}^4}\f{1}{b(x)}dx\ri)^{1/2}
  \le(\int_{\mathbb{R}^4}b(x)u^2dx\ri)^{1/2}\leq \le(\int_{\mathbb{R}^4}\f{1}{b(x)}dx\ri)^{1/2}\|u\|_E.$$
  For any $\gamma:1<\gamma<2$, there holds
  $$\int_{\mathbb{R}^4}|u|^\gamma dx\leq \int_{\mathbb{R}^4}(|u|+u^2)dx\leq \le(\int_{\mathbb{R}^4}\f{1}{b(x)}dx\ri)^{1/2}\|u\|_E
  +\f{1}{b_0}\|u\|_E^2,$$
  where $b_0$ is given by $(A_1)$.
  Thus we get continuous embedding $E\hookrightarrow
  L^q(\mathbb{R}^4)$ for all $q\geq 1$.

  To prove that the above embedding is also compact, take a sequence of functions
  $(u_k)\subset E$ such that $\|u_k\|_E\leq
  C$ for all $k$, we must prove that up to a subsequence there exists some $u\in E$ such that $u_k$
  convergent to  $u$ strongly in $L^q(\mathbb{R}^4)$
  for all $q\geq 1$. Without loss of generality, up to a subsequence, we may assume
  \be\label{bdd}\le\{\begin{array}{lll}
  u_k\rightharpoonup u\quad&{\rm weakly\,\,in}\quad E\\[1.5ex]
  u_k\rightarrow u\quad&{\rm strongly\,\,in}\quad L^q_{\rm loc}(\mathbb{R}^4),\,\,\forall q\geq 1\\[1.5ex]
  u_k\ra u\quad&{\rm almost\,\,everywhere\,\,in}\quad \mathbb{R}^4.
  \end{array}\ri.
  \ee
   In view of $(A_2)$, for any $\epsilon>0$, there exists $R>0$ such that
  $$\le(\int_{|x|> R}\f{1}{b}dx\ri)^{1/2}<\epsilon.$$
  Hence
  \be\label{res}\int_{|x|>R}|u_k-u|dx\leq \le(\int_{|x|> R}\f{1}{b}dx\ri)^{1/2}
  \le(\int_{|x|> R}b|u_k-u|dx\ri)^{1/2}\leq \epsilon\|u_k-u\|_E\leq
  C\epsilon.\ee
  On the other hand, it follows from (\ref{bdd}) that $u_k\ra u$ strongly in $L^1(\mathbb{B}_R)$.
  This together with (\ref{res}) gives
  $$\limsup_{k\ra \infty}\int_{\mathbb{R}^4}|u_k-u|dx\leq C\epsilon.$$
  Since $\epsilon$ is arbitrary, we obtain
  $$\lim_{k\ra \infty}\int_{\mathbb{R}^4}|u_k-u|dx=0.$$
  For $q>1$, it follows from the continuous embedding $E\hookrightarrow L^s(\mathbb{R}^4)$ ($s\geq 1$)
  that
  \bna
  \int_{\mathbb{R}^4}|u_k-u|^qdx&=&\int_{\mathbb{R}^4}|u_k-u|^\f{1}{2}|u_k-u|^{(q-\f{1}{2})}dx\\
  &\leq&\le(\int_{\mathbb{R}^4}|u_k-u|dx\ri)^{1/2}\le(\int_{\mathbb{R}^4}|u_k-u|^{2q-1}dx\ri)^{1/2}\\
  &\leq&C\le(\int_{\mathbb{R}^4}|u_k-u|dx\ri)^{1/2}\ra 0
  \ena
  as $k\ra\infty$. This concludes the lemma. $\hfill\Box$
  \\

  \noindent{\bf Lemma 3.7.} {\it Assume that $(A_1)$, $(A_2)$, $(H_1)$, $(H_2)$ and $(H_3)$ are satisfied.
  Let $(u_n)\subset E$ be an arbitrary Palais-Smale sequence of $J_\epsilon$, i.e.,
   \be\label{PS}J_{\epsilon}(u_n)\ra c,\,\, J^\prime_{\epsilon}(u_n)\ra 0 \,\, {\rm in}\,\,
   E^*
   \,\,{\rm as}\,\, n\ra \infty.\ee
   Then there exist a subsequence of $(u_n)$ (still denoted by $(u_n)$) and $u\in E$ such that
   $u_n\rightharpoonup u$ weakly in $E$, $u_n\ra u$ strongly in
   $L^q(\mathbb{R}^4)$ for all $q\geq 1$, and
   \bna\le\{\begin{array}{lll}
   \f{f(x,\,u_n)}{|x|^\beta}\ra \f{f(x,\,u)}{|x|^\beta}\,\,{\rm strongly\,\,in}\,\,
   L^1(\mathbb{R}^4)\\[1.5ex]
    \f{F(x,\,u_n)}{|x|^\beta}\ra \f{F(x,\,u)}{|x|^\beta}\,\, {\rm strongly\,\,in}\,\,
   L^1(\mathbb{R}^4).
   \end{array}\ri.\ena
   Furthermore $u$ is a weak solution of (\ref{problem}).\\
   }

   \noindent{\it Proof.} Assume $(u_n)$ is a Palais-Smale sequence of $J_\epsilon$. By $(\ref{PS})$, we have
   \bea\label{1}
   \f{1}{2}\int_{\mathbb{R}^4}\le(|\Delta u_n|^2+a|\nabla u_n|^2+bu_n^2\ri)dx-
   \int_{\mathbb{R}^4}\f{F(x,u_n)}{|x|^\beta}dx-\epsilon\int_{\mathbb{R}^4}hu_ndx\ra
   c\,\,{\rm as}\,\,n\ra\infty,\\\label{2}
   \le|\int_{\mathbb{R}^4}(\Delta u_n\Delta\psi+a\nabla
   u_n\nabla\psi+bu_n\psi)dx-\int_{\mathbb{R}^4}\f{f(x,u_n)}{|x|^\beta}\psi
   dx-\epsilon\int_{\mathbb{R}^4}h\psi dx\ri|\leq \tau_n\|\psi\|_E
   \eea
   for all $\psi\in E$, where $\tau_n\ra 0$ as $n\ra\infty$. By $(H_2)$, $0\leq \mu F(x,u_n)\leq
   u_nf(x,u_n)$ for some $\mu>2$. Taking $\psi=u_n$ in (\ref{2}) and multiplying
   (\ref{1}) by $\mu$, we have
   \bna
   \le(\f{\mu}{2}-1\ri)\|u_n\|_E^2&\leq&\le(\f{\mu}{2}-1\ri)\|u_n\|_E^2-\int_{\mathbb{R}^4}\f{\mu F(x,u_n)-f(x,u_n)u_n}{|x|^\beta}dx\\
   &\leq&\mu |c|+\tau_n\|u_n\|_E+(\mu+1)\epsilon\|h\|_{E^*}\|u_n\|_E
   \ena
   Therefore $\|u_n\|_E$ is bounded. It then follows from  (\ref{1}), (\ref{2}) that
   $$\int_{\mathbb{R}^4}\f{f(x,u_n)u_n}{|x|^\beta}
   dx\leq C,\quad \int_{\mathbb{R}^4}\f{F(x,u_n)}{|x|^\beta}dx\leq
   C.$$
   Notice that $f(x,u_n)u_n\geq 0$ and $F(x,u_n)\geq
   0$. By Lemma 3.6, up to a subsequence, $u_n\ra u$ strongly in $L^q(\mathbb{R}^4)$ for some $u\in E$, $\forall q\geq 1$.
    Then we obtain by applying Lemma 3.5 (here $N=4$ and $\phi=|x|^{-\beta}$),
   \be\label{f-conv}\lim_{n\ra\infty}\int_{\mathbb{R}^4}\f{|f(x,\,u_n)-f(x,\,u)|}{|x|^\beta}dx=0.\ee
   By $(H_1)$ and $(H_3)$, there exist constants $c_1$, $c_2>0$
   such that
   $$F(x,u_n)\leq c_1u_n^2+c_2|f(x,u_n)|.$$
   In view of (\ref{f-conv}) and Lemma 3.6, it follows from the generalized Lebesgue's dominated convergence
   theorem
   $$\lim_{n\ra\infty}\int_{\mathbb{R}^4}\f{|F(x,\,u_n)-F(x,\,u)|}{|x|^\beta}dx=0.$$
   Finally passing to the limit $n\ra \infty$ in $(\ref{2})$, we have
   $$\int_{\mathbb{R}^4}\le(\Delta u\Delta\psi+a\nabla
   u\nabla\psi+bu\psi\ri)dx-\int_{\mathbb{R}^4}\f{f(x,u)}{|x|^\beta}\psi
   dx-\epsilon\int_{\mathbb{R}^4}h\psi dx=0$$
   for all $\psi\in C_0^\infty(\mathbb{R}^4)$, which is dense in $E$. Hence $u$ is a weak solution of
   $(\ref{problem})$. $\hfill\Box$

 \subsection{Nontrivial solution\\}

 In this subsection, we will prove Theorem 1.3.
 It suffices to prove that the functional $J$ has a nontrivial
 critical point in the function space $E$.\\

 \noindent{\it Proof of Theorem 1.3.} Notice that $J_\epsilon$ becomes $J$ when $\epsilon=0$.
 By Lemma 3.1 and Lemma 3.2, $J$ satisfies all the hypothesis of
 the mountain-pass theorem except for the Palais-Smale condition: $J\in
   \mathcal{C}^1(E,\mathbb{R})$; $J(0)=0$; $J(u)\geq \delta>0$ when
  $\|u\|_E=r$; $J(e)<0$ for some  $e\in E$  with
   $\|e\|_E>r$.
    Then using the mountain-pass theorem
   without the Palais-Smale condition \cite{Rabin}, we can find a sequence
   $(u_n)$ of $E$ such that
   $$J(u_n)\ra c>0,\quad J^\prime(u_n)\ra 0\,\,{\rm in}\,\, E^*,$$
   where
   $$c=\min_{\gamma\in\Gamma}\max_{u\in\gamma}J(u)\geq \delta$$
   is the mountain-pass level of $J$, where $\Gamma=\{\gamma\in\mathcal{C}([0,1],E): \gamma(0)=0,
   \gamma(1)=e\}$. This is equivalent to saying
   \bea\label{01}
   \f{1}{2}\int_{\mathbb{R}^4}\le(|\Delta u_n|^2+a|\nabla u_n|^2+bu_n^2\ri)dx-
   \int_{\mathbb{R}^4}\f{F(x,u_n)}{|x|^\beta}dx\ra
   c\,\,{\rm as}\,\,n\ra\infty,\\\label{02}
   \le|\int_{\mathbb{R}^4}(\Delta u_n\Delta\psi+a\nabla
   u_n\nabla\psi+bu_n\psi)dx-\int_{\mathbb{R}^4}\f{f(x,u_n)}{|x|^\beta}\psi
   dx\ri|\leq \tau_n\|\psi\|_E
   \eea
   for $\psi\in E$, where $\tau_n\ra 0$ as $n\ra\infty$. By Lemma 3.7 with $\epsilon=0$,
   up to a subsequence, there holds
   \be\label{Fcon}\le\{\begin{array}{lll}
   u_n\rightharpoonup u \,\,{\rm weakly\,\, in}\,\, E\\[1.5ex]  u_n\ra
   u \,\,{\rm strongly\,\, in}\,\, L^q(\mathbb{R}^4),\,\,\forall q\geq
   1\\[1.5ex]
   \lim\limits_{n\ra\infty}\int_{\mathbb{R}^4}
   \f{F(x,u_n)}{|x|^\beta}dx=\int_{\mathbb{R}^4}\f{F(x,u)}{|x|^\beta}dx\\
   [1.5ex] u\,\,{\rm is\,\, a\,\,weak\,\,solution\,\,of}\,\,(\ref{h-zero}).\end{array}
   \ri.\ee

   Now suppose $u\equiv0$. Since $F(x,0)= 0$ for all $x\in \mathbb{R}^4$,  it follows from
   (\ref{01}) and (\ref{Fcon}) that
   \be\label{N}\lim_{n\ra\infty}\|u_n\|_E^2=2c>0.\ee
   Thanks to the hypothesis $(H_4)$, we have $0<c<16\pi^2/\alpha_0$ by applying Lemma 3.4.
   Thus there exists some $\epsilon_0>0$ and $N>0$ such that
   $\|u_n\|_E^2\leq 32\pi^2/\alpha_0-\epsilon_0$ for all $n>N$.
    Choose $q>1$ sufficiently close to $1$ such
   that $q\alpha_0\|u_n\|_E^2\leq 32\pi^2-\alpha_0\epsilon_0/2$ for all
   $n>N$. By $(H_1)$,
   $$|f(x,u_n)u_n|\leq b_1u_n^2+b_2|u_n|^{1+\gamma}\le(e^{\alpha_0u_n^2}-1\ri).$$
   It follows from the H\"older inequality, Lemma 2.1 and Theorem
   1.1 that
   \bna
   \int_{\mathbb{R}^4}\f{|f(x,u_n)u_n|}{|x|^\beta}dx&\leq&b_1\int_{\mathbb{R}^4}
   \f{u_n^2}{|x|^\beta}dx+b_2\int_{\mathbb{R}^4}\f{|u_n|^{1+\gamma}\le(e^{\alpha_0u_n^2}-1\ri)}{|x|^\beta}dx\\
   &\leq&b_1\int_{\mathbb{R}^4}\f{u_n^2}{|x|^\beta}dx+b_2\le(\int_{\mathbb{R}^4}\f{|u_n|^{(1+\gamma)q^\prime}}{|x|^\beta}dx\ri)^{1/{q^\prime}}
   \le(\int_{\mathbb{R}^4}\f{\le(e^{\alpha_0u_n^2}-1\ri)^q}{|x|^\beta}dx\ri)^{1/{q}}\\
   &\leq&b_1\int_{\mathbb{R}^4}\f{u_n^2}{|x|^\beta}dx+C\le(\int_{\mathbb{R}^4}\f{|u_n|^{(1+\gamma)q^\prime}}{|x|^\beta}dx\ri)^{1/{q^\prime}}
   \ra 0\quad{\rm as}\quad n\ra\infty.
   \ena
   Here we also used (\ref{Fcon}) (precisely $u_n\ra u$ in
   $L^s(\mathbb{R}^4)$ for all $s\geq 1$) in the last step of the above estimates.
    Inserting this into (\ref{02}) with $\psi=u_n$, we have
   $$\|u_n\|_E\ra 0\quad{\rm as}\quad n\ra\infty,$$
   which contradicts (\ref{N}). Therefore $u\not\equiv 0$ and we obtain a nontrivial
   weak solution of (\ref{h-zero}). $\hfill\Box$\\

   \subsection{Proof of Theorem 1.4\\}

   The proof of Theorem 1.4 is similar to the first part of that of Theorem 1.3.
   By Theorem 1.1, Lemma 3.1 and Lemma 3.2, there exists $\epsilon_1>0$ such that when $0<\epsilon<\epsilon_1$,
   $J_{\epsilon}$ satisfies all the
   hypothesis of the mountain-pass theorem except for the
   Palais-Smale condition: $J_{\epsilon}\in
   \mathcal{C}^1(E,\mathbb{R})$; $J_{\epsilon}(0)=0$; $J_{\epsilon}(u)\geq \vartheta_\epsilon>0$ when
  $\|u\|_E=r_\epsilon$; $J_{\epsilon}(e)<0$ for some  $e\in E$  with
   $\|e\|_E>r_\epsilon$.
    Then using the mountain-pass theorem
   without the Palais-Smale condition, we can find a sequence
   $(u_n)$ of $E$ such that
   $$J_{\epsilon}(u_n)\ra c>0,\quad J^\prime_{\epsilon}(u_n)\ra 0\,\,{\rm in}\,\, E^*,$$
   where
   $$c=\min_{\gamma\in\Gamma}\max_{u\in\gamma}J_{\epsilon}(u)\geq \vartheta_\epsilon$$
   is the mountain-pass level of $J_{\epsilon}$, where $\Gamma=\{\g\in\mathcal{C}([0,1],E): \g(0)=0,
   \g(1)=e\}$.
   By Lemma 3.7, there exists a subsequence of $(u_n)$ converges weakly  to a solution of
   (\ref{problem}) in $E$. $\hfill\Box$

  \subsection{A weak solution with negative energy\\}

  In this subsection, we will prove Theorem 1.5 by using the Ekeland's
   variational principle \cite{Mawin}. Let us first give
  two technical lemmas.\\

\noindent{\bf Lemma 3.8.} {\it Assume $(A_1)$ holds. If $(u_n)$ is a sequence in $E$ such
that
   \be\label{38}\lim_{n\ra\infty}\|u_n\|_E^2<\f{32\pi^2}{\alpha_0}\le(1-\f{\beta}{4}\ri),\ee
   then ${(e^{\alpha_0u_n^2}-1)}/{|x|^\beta}$ is bounded
   in $L^q(\mathbb{R}^4)$ for some $q>1$.}\\

   \noindent{\it Proof.} By $(A_1)$, we have
    $$\|u_n\|_E^2\geq\int_{\mathbb{R}^4}\le(|\Delta u_n|^2+a_0|\nabla u_n|^2+b_0u_n^2\ri)dx.$$
   Denote $v_n=u_n/\|u_n\|_E$. Then
    $$\int_{\mathbb{R}^4}\le(|\Delta v_n|^2+a_0|\nabla v_n|^2+b_0v_n^2\ri)dx\leq 1.$$
    By Theorem 1.1, for any $\alpha<32\pi^2(1-\beta/4)$, $0\leq\beta<4$, there holds
    \be\label{the}\int_{\mathbb{R}^4}\f{e^{\alpha v_n^2}-1 }{|x|^\beta}dx\leq
    C(\alpha,\beta).\ee
    By (\ref{38}), one can choose $q>1$ sufficiently close to $1$ such that
    \be\label{39}\lim_{n\ra\infty}\alpha_0q\|u_n\|_E^2< 32\pi^2(1-\beta
    q/4).\ee
    Combining Lemma 2.1, (\ref{the}) and (\ref{39}), we obtain
   \bna
   \int_{\mathbb{R}^4}\f{(e^{\alpha_0u_n^2}-1)^q}{|x|^{\beta
   q}}dx\leq\int_{\mathbb{R}^4}\f{e^{\alpha_0q\|u_n\|_E^2v_n^2}-1}{|x|^{\beta
   q}}dx\leq C
   \ena
   for some constant $C$. Thus ${(e^{\alpha_0u_n^2}-1)}/{|x|^\beta}$ is bounded
   in $L^q(\mathbb{R}^4)$. $\hfill\Box$\\

 \noindent{\bf Lemma 3.9.} {\it Assume $(A_1)$, $(A_2)$, $(H_1)$ are satisfied and $(u_n)$ is a Palais-Smale sequence for $J_\epsilon$ at any level
 with
 \be\label{assu}\liminf_{n\ra\infty}\|u_n\|_E^2<\f{32\pi^2}{\alpha_0}\le(1-\f{\beta}{4}\ri).\ee
 Then $(u_n)$ has a subsequence converging strongly to a solution of
 (\ref{problem}).}\\

 \noindent{\it Proof.} By (\ref{assu}), up to a subsequence, $(u_n)$ is bounded in $E$. In view of Lemma 3.6,
 without loss of generality we can assume
 \be\label{strong}\le\{\begin{array}{lll}
  u_n\rightharpoonup u_0\quad{\rm weakly\,\,in}\quad E\\[1.5ex]
  u_n\ra u_0\quad{\rm strongly\,\,in}\quad
  L^q(\mathbb{R}^4),\,\,\forall q\geq 1\\[1.5ex]
  u_n\ra u_0\quad{\rm almost\,\,everywhere\,\,in}\quad \mathbb{R}^4.
 \end{array}\ri.\ee
 Since $(u_n)$ is a Palais-Smale sequence for $J_\epsilon$, we have
 $J^\prime_\epsilon(u_n)\ra 0$ in $E^*$, particularly
 \bea{\nonumber}
 \int_{\mathbb{R}^4}\le(\Delta u_n\Delta(u_n-u_0)+a\nabla
 u_n\nabla(u_n-u_0)+bu_n(u_n-u_0)\ri)dx\\\label{p-s}-\int_{\mathbb{R}^4}
 \f{f(x,u_n)(u_n-u_0)}{|x|^\beta}dx-\epsilon\int_{\mathbb{R}^4}h(u_n-u_0)dx\ra
 0
 \eea
 as $n\ra\infty$. In view of (\ref{strong}), we have
 \be\label{we}\int_{\mathbb{R}^4}\le(\Delta u_0\Delta(u_n-u_0)+a\nabla
 u_0\nabla(u_n-u_0)+bu_0(u_n-u_0)\ri)dx\ra 0\quad{\rm as}\quad n\ra\infty.\ee
 Subtracting (\ref{we}) from (\ref{p-s}), we obtain
 \be\label{subt}\|u_n-u_0\|_E^2=\int_{\mathbb{R}^4}
 \f{f(x,u_n)(u_n-u_0)}{|x|^\beta}dx+\epsilon\int_{\mathbb{R}^4}h(u_n-u_0)dx+o(1).\ee
 In view of $(H_1)$ and (\ref{assu}), one can see from Lemma 3.8 that $f(x,u_n)$ is bounded
 in $L^q(\mathbb{R}^4)$ for some $q>1$ sufficiently close to $1$.
 It then follows from (\ref{strong}) and the H\"older inequality that
 \bna
 \int_{\mathbb{R}^4}\f{f(x,u_n)(u_n-u_0)}{|x|^\beta}dx+\epsilon\int_{\mathbb{R}^4}h(u_n-u_0)dx\ra
 0
 \ena
 as $n\ra\infty$.
 Inserting this into (\ref{subt}), we conclude
 $u_n\ra u_0$ strongly in $E$.
 Since $J_\epsilon\in \mathcal{C}^1(E,\mathbb{R})$,
  $u_0$ is a weak solution of
(\ref{problem}). $\hfill\Box$
\\

  \noindent {\it Proof of Theorem 1.5.}
    Let $r_\epsilon$ be as in Lemma 3.2, namely
    $J_{\epsilon}(u)>0$ for all
    $u:\|u\|_E=r_\epsilon$ with $r_\epsilon\ra 0$ as $\epsilon\ra 0$.
     One can choose
   $\epsilon_2: 0<\epsilon_2<\epsilon_1$ such that when $0<\epsilon<\epsilon_2$,
   \be\label{r-small}r_\epsilon<\f{32\pi^2}{\alpha_0}\le(1-\f{\beta}{4}\ri).\ee
   Lemma 3.8 together with $(H_1)$ and $(H_2)$ implies that $J_{\epsilon}(u)\geq -C$ for all $u\in\overline{B}_{r_\epsilon}
   =\{u\in E: \|u\|_E\leq
   r_\epsilon\}$, where $r_\epsilon$ is given by (\ref{r-small}).
      On the other hand thanks to Lemma 3.3, there holds $\inf_{\|u\|_E\leq r_\epsilon}J_{\epsilon}(u)<0$.
   Since $\overline{B}_{r_\epsilon}$ is a complete metric space with
   the metric given by the norm of $E$, convex and the functional
   $J_{\epsilon}$ is of class $\mathcal{C}^1$ and bounded
   below on $\overline{B}_{r_\epsilon}$, thanks to the Ekeland's
   variational principle, there exists some sequence $(u_n)\subset
   \overline{B}_{r_\epsilon}$ such that
   $$J_{\epsilon}(u_n)\ra c_0=\inf_{\|u\|_E\leq r_\epsilon}J_{\epsilon}(u),$$
   and $$J_{\epsilon}^\prime(u_n)\ra 0\quad {\rm in}\quad E^*$$
   as $n\ra\infty$.
   Observing that $\|u_n\|_E\leq
   r_\epsilon$, in view of (\ref{r-small}) and Lemma 3.9, we conclude that there exists a subsequence of $(u_n)$ which
   converges to a solution $u_0$ of (\ref{problem}) strongly in $E$. Therefore $J_{\epsilon}(u_0)=c_0<0$. $\hfill\Box$\\

 \subsection{Multiplicity results\\}

 In this subsection, we will show that two solutions obtained in
 Theorem 1.4 and Theorem 1.5 are distinct under some assumptions, i.e., Theorem 1.6 holds.
 We need the following technical lemma:\\

 \noindent{\bf Lemma 3.10.} {\it Let $(w_n)$ be a sequence in $E$. Suppose $\|w_n\|_E=1$ and
 $w_n\rightharpoonup w_0$ weakly in $E$. Then for any $0<p<\f{1}{1-\|w_0\|_E^2}$
 \be\label{3.49}\sup_n\int_{\mathbb{R}^4}\f{e^{32\pi^2(1-\beta/4)pw_n^2}-1}{|x|^\beta}dx<\infty.\ee
 }
 \noindent{\it Proof.} Since $w_n\rightharpoonup w_0$ weakly in $E$ and $\|w_n\|_E=1$,
 we have
 \bea{\nonumber}
 \|w_n-w_0\|_E^2&=&\int_{\mathbb{R}^4}\le(|\Delta (w_n-w_0)|^2+a(x)|\nabla
 (w_n-w_0)|^2+b(x)(w_n-w_0)^2\ri)dx\\{\nonumber}
 &=&1+\|w_0\|_E^2-2\int_{\mathbb{R}^4}\le(\Delta w_n\Delta w_0+a(x)\nabla
 w_n\nabla w_0+b(x)w_nw_0\ri)dx\\{\label{lim}}
 &\ra& 1-\|w_0\|_E^2\quad{\rm as}\quad n\ra\infty.
 \eea
 If $w_0\equiv 0$, then (\ref{3.49}) is a consequence of Theorem 1.1. If
 $w_0\not\equiv 0$, using the H\"older inequality, Lemma 2.1, Theorem 1.1 and the
 inequality
 \bna
  rs-1\leq \f{1}{\mu}\le(r^\mu-1\ri)+\f{1}{\nu}\le(s^\nu-1\ri),
 \ena
 where $r\geq 0, s\geq 0,\mu>1,\nu>1,
  {1}/{\mu}+{1}/{\nu}=1$, we estimate
  \bea{\nonumber}
  \int_{\mathbb{R}^4}\f{e^{32\pi^2(1-\beta/4)pw_n^2}-1}{|x|^\beta}dx&\leq&
  \int_{\mathbb{R}^4}\f{e^{32\pi^2(1-\beta/4)p\le((1+\epsilon)(w_n-w_0)^2+
  (1+\epsilon^{-1})w_0^2\ri)}-1}{|x|^\beta}dx\\{\nonumber}&\leq&
  \f{1}{q}\int_{\mathbb{R}^4}\f{e^{32\pi^2(1-\beta/4)qp(1+\epsilon)(w_n-w_0)^2}-1}
  {|x|^\beta}dx\\\label{33}
  &&+\f{1}{q^\prime}\int_{\mathbb{R}^4}\f{e^{32\pi^2(1-\beta/4)q^\prime p(1+\epsilon^{-1})w_0^2}-1}
  {|x|^\beta}dx,
  \eea
  where $1/q+1/q^\prime=1$. Assume $0<p<{1}/{(1-\|w_0\|_E^2)}$. By (\ref{lim}),
  we can choose $q$ sufficiently close to $1$ and $\epsilon>0$ sufficiently
  small such that
  $$qp(1+\epsilon)\|w_n-w_0\|_E^2<1$$
  for large $n$.
  Recall that $\|u\|_E^2\geq \int_{\mathbb{R}^4}(|\Delta u|^2+a_0|\nabla
  u|^2+b_0u^2)dx$. Applying Theorem 1.1, we conclude the lemma from (\ref{33}). $\hfill\Box$\\

  We remark that similar results were obtained in \cite{Lu-Yang} for
  bi-Laplacian on bounded smooth domain $\Omega\subset\mathbb{R}^4$
  and in \cite{do-de} for Laplacian on the whole $\mathbb{R}^2$.\\

 \noindent{\it Proof of Theorem 1.6.} According to Theorem 1.4 and Theorem 1.5, under the assumptions of Theorem
 1.6, there exist
 sequences $(v_n)$ and $(u_n)$ in $E$ such that as $n\ra\infty$,
 \bea{\label{42}}
  &&v_n\rightharpoonup u_M\,\,{\rm weakly\,\,in}\,\,E,\quad J_\epsilon(v_n)\ra c_M>0,\quad
  |\langle
  J_\epsilon^\prime(v_n),\phi\rangle|\leq \gamma_n\|\phi\|_E\\
  &&{\label{41}}
  u_n\ra u_0\,\,{\rm strongly\,\,in}\,\,E,\quad J_\epsilon(u_n)\ra c_0<0,\quad |\langle
  J_\epsilon^\prime(u_n),\phi\rangle|\leq \tau_n\|\phi\|_E
 \eea
 with $\gamma_n\ra 0$ and $\tau_n\ra 0$, both $u_M$ and $u_0$ are nonzero weak solutions to (\ref{problem})
 since $h\not\equiv 0$ and $\epsilon>0$.
 Suppose $u_M=u_0$. Then $v_n\rightharpoonup u_0$ weakly in $E$ and
 thus
  $$\int_{\mathbb{R}^4}\le(\Delta v_n\Delta u_0+a\nabla v_n\nabla u_0+bv_nu_0\ri)dx\ra
  \|u_0\|_E^2$$
  as $n\ra\infty$. Using the H\"older inequality, we obtain
 $$\limsup_{n\ra\infty}\|v_n\|_E\geq \|u_0\|_E>0.$$
 On one hand, by Lemma 3.7, we have
 \be\label{cv}\int_{\mathbb{R}^4}\f{F(x,v_n)}{|x|^\beta}dx\ra\int_{\mathbb{R}^4}\f{F(x,u_0)}{|x|^\beta}dx
 \quad{\rm as}\quad n\ra\infty.\ee
 Here and in the sequel, we do not distinguish sequence and subsequence. On the other hand,
 it follows from Theorem 1.4 that $\|v_n\|_E$ is bounded. In view of Lemma 3.6, it holds
 \be\label{ch}
 \int_{\mathbb{R}^4}hv_ndx\ra \int_{\mathbb{R}^4}hu_0dx.\quad{\rm as}\quad n\ra\infty.
 \ee
 Inserting (\ref{cv}) and (\ref{ch}) into (\ref{42}), we obtain
 \be\label{cM}
 \f{1}{2}\|v_n\|_E^2=c_M+\int_{\mathbb{R}^4}\f{F(x,u_0)}{|x|^\beta}dx+\epsilon\int_{\mathbb{R}^4}hu_0dx+o(1),
 \ee
 where $o(1)\ra 0$ as $n\ra\infty$.
 In the same way, one can derive
 \be\label{c0}
 \f{1}{2}\|u_n\|_E^2=c_0+\int_{\mathbb{R}^4}\f{F(x,u_0)}{|x|^\beta}dx+\epsilon\int_{\mathbb{R}^4}hu_0dx+o(1).
 \ee
 Combining (\ref{cM}) and (\ref{c0}), we have
 \be\label{Lp}
 \|v_n\|_E^2-\|u_0\|_E^2=2\le(c_M-c_0+o(1)\ri).
 \ee
 Now we need to estimate $c_M-c_0$. By Lemma 3.4, there holds for sufficiently small $\epsilon>0$,
 $$\max_{t\geq 0}J_\epsilon(t\phi_n)<\f{16\pi^2}{\alpha_0}\le(1-\f{\beta}{4}\ri).$$
 Since $c_M$ is the mountain-pass level of $J_\epsilon$, we have
 $$c_M<\f{16\pi^2}{\alpha_0}\le(1-\f{\beta}{4}\ri).$$
 From the proof of Lemma 3.3, we know that $c_0\ra 0$ as $\epsilon\ra 0$ ($c_0$ depends on $\epsilon$).
 Noting that $c_M>0$ and $c_0<0$, we obtain for sufficiently small $\epsilon>0$,
 \be\label{c-c}0<c_M-c_0<\f{16\pi^2}{\alpha_0}\le(1-\f{\beta}{4}\ri).\ee
  Write
 $$w_n=\f{v_n}{\|v_n\|_E},\quad w_0=\f{u_0}{\le(\|u_0\|_E^2+2(c_M-c_0)\ri)^{1/2}}.$$
 It follows from (\ref{Lp}) and $v_n\rightharpoonup u_0$ weakly in
 $E$ that
 $w_n\rightharpoonup w_0$ weakly in $E$. Notice that
 $$\int_{\mathbb{R}^4}\f{e^{\alpha_0v_n^2}-1}{|x|^\beta}dx=\int_{\mathbb{R}^4}\f{e^{\alpha_0\|v_n\|_E^2
 w_n^2}-1}{|x|^\beta}dx.$$
 By (\ref{Lp}) and (\ref{c-c}), a straightforward calculation shows
 $$\lim_{n\ra\infty}\alpha_0\|v_n\|_E^2(1-\|w_0\|_E^2)<
 32\pi^2\le(1-\f{\beta}{4}\ri).$$
 Whence Lemma 3.10 implies that
 $e^{\alpha_0v_n^2}$ is bounded in $L^q(\mathbb{R}^4)$ for some
 $q>1$. By $(H_1)$,
 $$|f(x,v_n)|\leq b_1|v_n|+b_2|v_n|^\gamma\le(e^{\alpha_0v_n^2}-1\ri).$$
 Then the H\"older inequality and the continuous embedding $E\hookrightarrow L^p(\mathbb{R}^4)$
 for all $p\geq 1$ imply that
 $f(x,v_n)/|x|^\beta$ is bounded in $L^{q_1}(\mathbb{R}^4)$ for
 some $q_1$: $1<q_1<q$.
  This together with Lemma
 3.6 and the H\"older inequality gives
 \be\label{f-c}
 \le|\int_{\mathbb{R}^4}\f{f(x,v_n)(v_n-u_0)}{|x|^\beta}dx\ri|\leq
 \le\|\f{f(x,v_n)}{|x|^\beta}\ri\|_{L^{q_1}(\mathbb{R}^4)}\le\|v_n-u_0\ri\|_{L^{{q_1^\prime}}(\mathbb{R}^4)}
 \ra 0,
 \ee
 where $1/q_1+1/q_1^\prime=1$.

 Taking $\phi=v_n-u_0$ in (\ref{42}), we have by using (\ref{f-c}) and
 Lemma 3.6 that
 \be\label{cc}
  \int_{\mathbb{R}^4}\le(\Delta v_n\Delta(v_n-u_0)+a\nabla
  v_n\nabla(v_n-u_0)+bv_n(v_n-u_0)\ri)dx
  \ra 0.
 \ee
 However the fact $v_n\rightharpoonup u_0$ weakly in $E$ implies
 \be\label{cc1}
  \int_{\mathbb{R}^4}\le(\Delta u_0\Delta(v_n-u_0)+a\nabla
  u_0\nabla(v_n-u_0)+bu_0(v_n-u_0)\ri)dx
  \ra 0.
 \ee
 Subtracting (\ref{cc1}) from (\ref{cc}), we have $\|v_n-u_0\|_E^2\ra
 0$. Since $J_\epsilon\in\mathcal{C}^1(E,\mathbb{R})$, we have
 $$J_\epsilon(v_n)\ra J_\epsilon(u_0)=c_0,$$
 which contradicts $J_\epsilon(v_n)\ra c_M>c_0$. This completes the proof of Theorem 1.6.
 $\hfill\Box$

  \end{document}